\documentclass[11pt]{article}
\usepackage{epic,latexsym,amssymb}
\usepackage{color}
\usepackage{tikz}
\usepackage{amsfonts,epsf,amsmath}

\newtheorem{thm}{Theorem}
\newtheorem{lem}{Lemma}

\newtheorem{defn}{Definition}
\newtheorem{claim}{Claim}

\newcommand{\qed}{$\Box$}

\newcommand{\smallqed}{{\tiny ($\Box$)}}

\newcommand{\ii}{\iota}

\newcommand{\proof}{\noindent\textbf{Proof. }}

\let\oldenumerate\enumerate
\renewcommand{\enumerate}{
  \oldenumerate
  \setlength{\itemsep}{0pt}
  \setlength{\parskip}{0pt}
  \setlength{\parsep}{0pt}
}

\begin{document}

\title{Extensions of the Art Gallery Theorem}

\author{Peter Borg\\[5mm]
\normalsize Department of Mathematics \\
\normalsize Faculty of Science \\
\normalsize University of Malta\\
\normalsize Malta\\
\normalsize \texttt{peter.borg@um.edu.mt} 
\and 
Pawaton Kaemawichanurat \\ [5mm]
\normalsize Mathematics and Statistics with Applications (MaSA) \\
\normalsize and Department of Mathematics, Faculty of Science\\
\normalsize King Mongkut's University of Technology Thonburi \\
\normalsize Bangkok, Thailand\\
\normalsize \texttt{pawaton.kae@kmutt.ac.th}
}


\date{}
\maketitle

\begin{abstract}
Several domination results have been obtained for maximal outerplanar graphs (mops). The classical domination problem is to minimize the size of a set $S$ of vertices of an $n$-vertex graph $G$ such that $G - N[S]$, the graph obtained by deleting the closed neighborhood of $S$, contains no vertices. In the proof of the Art Gallery Theorem, Chv\'{a}tal showed that the minimum size, called the domination number of $G$ and denoted by $\gamma(G)$, is at most $n/3$ if $G$ is a mop. Here we consider a modification by allowing $G - N[S]$ to have a maximum degree of at most $k$. Let $\ii_k(G)$ denote the size of a smallest set $S$ for which this is achieved.
%
%
%
If $n \le 2k+3$, then trivially $\ii_k(G) \leq 1$. Let $G$ be a mop on $n \ge \max\{5,2k+3\}$ vertices, $n_2$ of which are of degree~$2$. Upper bounds on $\ii_k(G)$ have been obtained for $k = 0$ and $k = 1$, namely $\ii_{0}(G) \le \min\{\frac{n}{4},\frac{n+n_2}{5},\frac{n-n_2}{3}\}$ and $\ii_1(G) \le \min\{\frac{n}{5},\frac{n+n_2}{6},\frac{n-n_2}{3}\}$.
We prove that $\ii_{k}(G) \le \min\{\frac{n}{k+4},\frac{n+n_2}{k+5},\frac{n-n_2}{k+2}\}$ for any $k \ge 0$. For the original setting of the Art Gallery Theorem, the argument presented yields that if 
an art gallery has exactly $n$ corners and
at least one of every $k + 2$ consecutive corners must be visible to at least one guard, then the number of guards needed is at most $n/(k+4)$. We also prove that $\gamma(G) \le \frac{n - n_2}{2}$ unless $n = 2n_2$, $n_2$ is odd, and $\gamma(G) = \frac{n - n_2 + 1}{2}$. 
Together with the inequality $\gamma(G) \le \frac{n+n_2}{4}$, obtained by Campos and Wakabayashi and independently by Tokunaga, this improves Chv\'{a}tal's bound. The bounds are sharp.
\end{abstract}


\section{Introduction}
Unless stated otherwise, we use small letters such as $x$ to denote integers or elements of a set, and capital letters such as $X$ to denote sets or graphs. The set of positive integers is denoted by $\mathbb{N}$. For $m, n \in \{0\} \cup \mathbb{N}$, the set $\{i \in \mathbb{N} \colon m \leq i \leq n\}$ is denoted by $[m,n]$. We abbreviate $[1,n]$ to $[n]$. Note that $[0]$ is the empty set $\emptyset$. For a set $X$, the set of $k$-element subsets of $X$ is denoted by ${X \choose k}$. 
Arbitrary sets are assumed to be finite.

Let $G$ be a simple graph with vertex set $V(G)$ and edge set $E(G)$. We may represent an edge $\{v,w\}$ by $vw$. The \emph{order of $G$} is $|V(G)|$.  
We say that $G$ is an \emph{$n$-vertex graph} if its order is $n$. The \emph{open neighborhood} $N_G(v)$ of a vertex $v$ of $G$ is the set of neighbors of $v$, that is, $N_G(v) = \{w \in V(G) \colon vw \in E(G)\}$. The \emph{degree} $d_G(v)$ of $v$ is $|N_G(v)|$. The \emph{maximum degree} $\Delta(G)$ of $G$ is $\max\{d_G(v) \colon v \in V(G)\}$. The \emph{closed neighborhood} $N_G[v]$ of $v$ is the set $N_G(v) \cup \{v\}$.
For $U \subseteq V(G)$, the \emph{closed neighborhood} $N_G[U]$ of $U$ is $\bigcup_{v \in U} N_G[v]$. The subgraph of $G$ \emph{induced by $U$} is denoted by $G[U]$, that is, $V(G[U]) = U$ and $E(G[U]) = E(G) \cap {U \choose 2}$. The subgraph of $G$ obtained by deleting the vertices in $U$ from $G$ is denoted by $G - U$, that is, $G-U = G[V(G) \setminus U]$. We may abbreviate $G - \{v\}$ to $G - v$. 
We use the standard notation $K_k$, $P_k$, $C_k$, and $K_{1,k}$ for the $k$-vertex \emph{complete graph}, 
the $k$-vertex \emph{path}, 
the $k$-vertex \emph{cycle}, 
and the $(k+1)$-vertex \emph{star} ($E(K_{1,k}) = \{uv \colon v \in V(K_{1,k}) \setminus \{u\}\}$ for some $u \in V(K_{1,k})$), respectively. If $G$ is a $k$-vertex path and $E(G) = \{v_iv_{i+1} \colon i \in [k-1]\}$, then we represent $G$ by $v_1v_2 \dots v_k$. If $G$ is a $k$-vertex cycle and $E(G) = \{v_iv_{i+1} \colon i \in [k-1]\} \cup \{v_kv_1\}$, then we represent $G$ by $v_1v_2 \dots v_kv_1$.

A subset $S$ of $V(G)$ is a \emph{dominating set of $G$} if each vertex in $V(G) \setminus S$ is adjacent to at least one vertex in $S$ (that is, $N_G[S] = V(G)$). The classical domination problem is to minimize the size of a dominating set (see \cite{C, CH, HHS, HHS2, HL, HL2}). Caro and Hansberg~\cite{CaHa17} introduced an appealing generalization that has also been treated in \cite{Borg_conn, Borg1, BFK, BFK2, BK, FK, KaJi, Y, YW, ZW}. They relaxed the domination condition by considering a subset $S$ of $V(G)$ such that $G - N_{G}[S]$ contains no forbidden subgraph. Given a graph $H$, $S$ is called an \emph{$H$-isolating set of $G$} if $G - N_{G}[S]$ does not contain a copy of $H$. The \emph{$H$-isolation number of $G$} is the size of a smallest $H$-isolating set of $G$ and is denoted by $\ii(G, H)$. The \emph{domination number of $G$} is the size of a smallest dominating set of $G$ and is denoted by $\gamma(G)$. Note that $S$ is a dominating set if and only if it is a $K_{1}$-isolating set; thus, $\gamma(G) = \ii(G, K_1)$. We are interested in the case $H = K_{1,k+1}$. Note that, for $k \geq 0$, $S$ is a $K_{1,k+1}$-isolating set of $G$ if and only if $\Delta(G-N_G[S]) \leq k$. We abbreviate $\ii(G, K_{1,k+1})$ to $\ii_{k}(G)$. Since $K_{1,0}$ is the graph $K_1$, we have $\gamma(G) = \ii_{-1}(G)$.

A \emph{triangulated disc} is a plane graph whose interior faces are triangles and whose exterior face (the unbounded face) is bounded by a simple cycle. A \emph{maximal outerplanar graph}, or a \emph{mop}, is a triangulated disc $G$ such that the boundary of the exterior face of $G$ contains all the vertices of $G$. 
O'Rourke~\cite{Ro87} pointed out that every mop has a unique Hamiltonian cycle. Thus, the Hamiltonian cycle of a mop is the boundary of the mop. This paper's notation and terminology on mops follows that of~\cite{LeZuZy17}; in particular, each edge of the Hamiltonian cycle of a mop is called a \emph{Hamiltonian edge}, while any other edge of the mop is called a \emph{diagonal}. For $n \geq 3$, the \emph{fan} $F_n$ is the mop obtained from $P_{n-1}$ by adding a new vertex $v$ and joining it to every vertex of $P_{n-1}$. We say that $v$ is the \emph{center of $F_n$}.

Domination in mops has been extensively studied since 1975. Chv\'{a}tal's proof of his classical result referred to as the Art Gallery Theorem (AGT) \cite{Ch75} established that the domination number of any $n$-vertex mop is at most $n/3$. This also follows from Fisk's elegant proof \cite{Fi78} (included in \cite{AZ}) of AGT, and was proved directly by Matheson and Tarjan~\cite{MaTa96}. For results on other types of domination in mops, we refer the reader to \cite{CaCaHeMa16,DoHaJo16,DoHaJo17,HeKa18}. Caro and Hansberg~\cite{CaHa17} proved that the $K_{1, 1}$-isolation number of a mop of order $n \ge 4$ is at most $n/4$. Borg and Kaemawichanurat~\cite{BK} proved that the $K_{1, 2}$-isolation number of a mop of order $n \ge 5$ is at most $n/5$.

\parskip=0pt

\begin{thm}
\label{t:known1} \setlength{\baselineskip}{0.6cm}
Let $G$ be an $n$-vertex mop. \\
(a) {\rm (\cite{Ch75})} If $n \ge 3$, then $\gamma(G) \le \frac{n}{3}$. \\
(b) {\rm (\cite{CaHa17})} If $n \ge 4$, then $\ii_0(G) \le \frac{n}{4}$. \\
(c) {\rm (\cite{BK})} If $n \ge 5$, then $\ii_1(G) \le \frac{n}{5}$.\\
Moreover, the bounds are sharp.
\end{thm}
%
When we say that a bound is sharp, we mean that there are infinitely many values of $n$ for which the bound is attained. For each of the bounds in Theorem~\ref{t:known1}, the bound is attained for each $n$ such that the bound is an integer.

\parskip=4pt
The following sharp upper bounds in terms of the order and the number of vertices of degree $2$ have also been established.
\parskip=0pt

\begin{thm}
\label{t:known2.0} \setlength{\baselineskip}{0.6cm}
If $G$ is a mop of order $n \geq 3$ and has exactly $n_2$ vertices of degree $2$, then  \\
(a) {\rm (\cite{CaWa13,To13})} 
$\gamma(G)  \le \frac{n + n_2}{4}$,\\
(b) {\rm (\cite{KaJi})} 
$\ii_{0}(G)  \le \frac{n + n_2}{5}$,\\
(c) {\rm (\cite{BK})} 
$\ii_{1}(G)  \le \frac{n + n_2}{6}$.\\
Moreover, the bounds are sharp.
\end{thm}

\begin{thm}
\label{t:known2}
\setlength{\baselineskip}{0.6cm}
If $G$ is a mop of order $n \geq 5$ and has exactly $n_2$ vertices of degree $2$, then  \\
(a) {\rm (\cite{KaJi})} 
$\ii_{0}(G)  \le \frac{n - n_2}{3}$,\\
(b) {\rm (\cite{BK})} 
$\ii_{1}(G)  \le \frac{n - n_2}{3}$.\\
Moreover, the bounds are sharp.
\end{thm}

\section{Main results}\label{main}
In this paper, we mostly establish sharp upper bounds on the $K_{1, k + 1}$-isolation number of a mop in terms of its order and the number of vertices of degree~$2$ for any $k \ge -1$. Our new results are presented in this section and in Section~\ref{AGTextsection}, and proved in Sections~\ref{S:main1} and \ref{AGTextsection}, respectively. The first three results have Theorems~\ref{t:known1} and \ref{t:known2.0} for $0 \leq k \leq 1$, and Theorem~\ref{t:known2}(b), as special cases. 

\begin{thm}
\label{t:main1}
If $k \geq 0$ and $G$ is a mop of order $n \ge k+4$, then
\[\ii_{k}(G) \leq \frac{n}{k + 4}.\]
Moreover, the bound is sharp.
\end{thm}
%
In Section~\ref{AGTextsection}, we present bounds that follow by the same argument used in the proof of Theorem~\ref{t:main1}. We address a relaxation of the condition in the Art Gallery Theorem that the whole polygon needs to be guarded. We show that if $k \ge -1$, a polygon has exactly $n \ge k+4$ corners, and at least one of every $k + 2$ consecutive corners must be visible to at least one guard, then the number of guards needed is at most $n/(k + 4)$ (see Theorem~\ref{AGTextthm}).

\begin{thm}
\label{t:main2}
If $k \geq 0$, $G$ is a mop of order $n \geq k+3$, and $n_2$ is the number of vertices of $G$ of degree~$2$, then
\[ \ii_{k}(G)  \le \frac{n + n_2}{k + 5}.\]
Moreover, the bound is sharp.
\end{thm}

\begin{thm}
\label{t:main3}
If $k \geq 1$, $G$ is a mop of order $n \geq 2k+3$, and $n_2$ is the number of vertices of $G$ of degree~$2$, then
\[ \ii_{k}(G)  \le \frac{n - n_2}{k + 2}.\]
Moreover, the bound is sharp.
\end{thm}

Note that, surprisingly, the sharp bound for $k = 0$ given by Theorem~\ref{t:known2} is not of the general form for $k \geq 1$ given by Theorem~\ref{t:main3}. 
%
Unlike Theorems~\ref{t:known1} and \ref{t:known2.0}, Theorem~\ref{t:known2} provides no bound for $\gamma(G)$ similar to those in its parts (a) and (b). The missing bound is provided by our next theorem. 

\begin{defn} \label{extreme} \emph{If $t \geq 3$, $x_1x_2 \dots x_{t}x_1$ is the unique Hamiltonian cycle $C$ of a mop $M$ contained by a mop $G$, $y_1, \dots, y_t$ are distinct vertices of $G$ in the exterior face of $M$, $N_G(y_i) = \{x_i, x_{i+1}\}$ for each $i \in [t-1]$, $N_G(y_{t}) = \{x_{t}, x_{1}\}$, $V(G) = V(M) \cup \{y_i \colon i \in [t]\}$, and $E(G) = E(M) \cup \bigcup_{i=1}^{t} \{y_ix \colon x \in N_G(y_i)\}$, then we call $G$ a \emph{$t$-extreme} mop or simply an \emph{extreme} mop.}
\end{defn}

Clearly, a mop $G$ as in Definition~\ref{extreme} can be constructed for any $t \ge 3$, and its Hamiltonian cycle is $x_1y_1x_2y_2 \dots x_ty_tx_1$. 

\parskip=4pt
For positive integers $a$ and $b$, let
\[ {\bf 1}(a,b)  = \left\{\begin{array}{ll}
                     1  \quad  & \mbox{if $b$ is odd and $a = 2b$},\\ [1mm]
                     0  \quad  & \mbox{otherwise}. 
               \end{array}\right. \]

\begin{thm}
\label{t:main4}
If $G$ is a mop of order $n \geq 4$ and has exactly $n_2$ vertices of degree $2$, then
\[\gamma(G)  \le \frac{n - n_2 + {\bf 1}(n,n_2)}{2}.\]
Moreover, 
equality holds if $G$ is extreme or $n = 4$.
\end{thm}
\parskip=0pt
Therefore, unlike Theorem~\ref{t:main3} for $k \ge 1$, Theorems~\ref{t:known2}(a) and \ref{t:main4} yield $\ii_k(G) \le \frac{n - n_2 - k {\bf 1}(n,n_2)}{k + 3}$ for $-1 \le k \le 0$ (recall that $\gamma(G) = \ii_{-1}(G)$) and $n \ge k+5$.

\parskip=4pt
Theorem~\ref{t:main4} enables us to improve the classical bound of Chv\'{a}tal in Theorem~\ref{t:known1}(a), using Theorem~\ref{t:known2.0}(a). The improved bound is given in Theorem~\ref{t:MAIN}, which also specifies a necessary condition and a sufficient condition for an $n$-vertex mop $G$ to attain Chv\'{a}tal's upper bound $\frac{n}{3}$ on $\gamma(G)$.

\parskip=0pt

\begin{defn} \label{special} \emph{If $t \geq 2$, $x_1x_2 \dots x_{2t}x_1$ is the unique Hamiltonian cycle $C$ of a mop $M$ contained by a mop $G$, $y_1, \dots, y_t$ are distinct vertices of $G$ in the exterior face of $M$, $N_G(y_1), \dots, N_G(y_t)$ are distinct edges of $C$, $x_i \in N_G[\{y_1, \dots, y_t\}]$ for each $i \in [2t]$ with $d_M(x_i) = 2$, $V(G) = V(M) \cup \{y_i \colon i \in [t]\}$, and $E(G) = E(M) \cup \bigcup_{i=1}^t \{y_ix \colon x \in N_G(y_i)\}$, then we call $G$ a \emph{$t$-special} mop or simply a \emph{special} mop. If $N_G(y_1), \dots, N_G(y_t)$ partition $V(C)$ (that is, their union is $V(C)$ and no two of them intersect), then we call $G$ an \emph{extra $t$-special} mop or simply an \emph{extra special} mop.}
\end{defn}

Clearly, an extra $t$-special mop can be constructed for any $t \ge 2$.

\begin{thm}
\label{t:MAIN}
If $G$ is a mop of order $n \geq 4$ and has exactly $n_2$ vertices of degree $2$, then \\
\[ \gamma(G)  \leq \left\{\begin{array}{ll}
                     \frac{n + n_{2}}{4} < \frac{n}{3}  \quad  & \mbox{if}~~n_{2} < \frac{n}{3},\\\\
                     \frac{n + n_2}{4} = \frac{n - n_2 + {\bf 1}(n,n_2)}{2} = \frac{n}{3}  \quad  & \mbox{if}~~n_{2} = \frac{n}{3},\\\\
                     \frac{n - n_2 + {\bf 1}(n,n_2)}{2} = \frac{n}{3}  \quad  & \mbox{if}~~(n,n_{2}) = (6,3),\\\\
                     \frac{n - n_2 + {\bf 1}(n,n_2)}{2} < \frac{n}{3}  \quad  & \mbox{if}~~n_2 > \frac{n}{3}~~\mbox{and}~~(n,n_2) \neq (6,3).
               \end{array}\right. \] \\
Moreover, the following assertions hold: \\
(a) The bound is sharp. \\
(b) If $\gamma(G) = \frac{n}{3}$, then $G$ is $n_2$-special or $3$-extreme. \\
(c) If $G$ is extra $n_2$-special, then $\gamma(G) = \frac{n}{3}$.
\end{thm}

Clearly, $|V(K_{1,k+1})| = k+2$. It immediately follows that, for any $n$-vertex mop $G$, $\ii_k(G) = 0$ if $n \leq k+1$, and $\ii_k(G) \leq 1$ if $n \leq \max\{k+4, 2k+3\}$. 
Theorems~\ref{t:known1}--\ref{t:main4} immediately give us the following complete solution for any other value of $n$.

\begin{thm}
\label{t:cor}
Let $G$ be an $n$-vertex mop, and let $n_2$ be the number of vertices of $G$ of degree~$2$. \\
(a) If $k \geq 1$ and $n \ge 2k+3$, then 
\[ \ii_{k}(G)  \leq \left\{\begin{array}{ll}
                     \frac{n + n_{2}}{k + 5}  \quad  & \mbox{if}~~n_{2} \leq \frac{n}{k + 4},\\\\
                     \frac{n}{k + 4}  \quad  & \mbox{if}~~\frac{n}{k + 4} \leq n_{2} \leq \frac{2n}{k + 4},\\\\
                     \frac{n - n_{2}}{k + 2}  \quad  & \mbox{if}~~n_2 \geq \frac{2n}{k + 4}.
               \end{array}\right. \]
(b) If $-1 \le k \le 0$ and $n \geq k+5$, then
\[ \ii_{k}(G)  \leq \left\{\begin{array}{ll}
                     \frac{n + n_{2}}{k+5}  \quad  & \mbox{if}~~n_{2} \leq \frac{n}{k+4},\\\\
                     \frac{n - n_{2} - k{\bf 1}(n,n_2)}{k+3}  \quad  & \mbox{if}~~n_{2} \geq \frac{n}{k+4}. 
               \end{array}\right. \]
%
%
%
Moreover, the bound is sharp.
\end{thm}




\section{Extremal constructions for Theorems~\ref{t:main1}--\ref{t:main3}} \label{extremalsection}

We now show that the bounds in Theorems~\ref{t:main1}--\ref{t:main3} are attainable. Theorems~\ref{t:known1}--\ref{t:known2} already establish this for $k \le 1$, so we settle $k \geq 2$.

\parskip=4pt
When $n/(k+4) < n_2 < 2n/(k+4)$, we have $n/(k+4) < \min\{(n+n_2)/(k+5), (n-n_2)/(k+2)\}$, that is, the bound in Theorem~\ref{t:main1} is better than those in Theorems~\ref{t:main2} and \ref{t:main3} for this range. Therefore, we will first show that the bound $n/(k+4)$ is attained in cases where $n/(k+4) < n_2 < 2n/(k+4)$. For an integer $t \ge 1$, let $F^1_{k + 4}, F^2_{k + 4}, \dots, F^{2t}_{k + 4}$ be $2t$ vertex-disjoint fans of order $k+4$. For $i \in [2t]$, let $x^{i}_{0}, x^{i}_{1}, \dots, x^{i}_{k+3}$ be the vertices of $F^i_{k + 4}$ with $x^{i}_{0}$ being the center and with $x^{i}_{1}$ and $x^{i}_{k+3}$ being the vertices of degree $2$. We extend the union of $F^1_{k + 4}, F^2_{k + 4}, \dots, F^{2t}_{k + 4}$ to a mop $A_{k,t}$ by adding edges on the $4t$ vertices $x^{1}_{1}, x^{1}_{2}, \dots, x^{t}_{1}, x^{t}_{2}, x^{t+1}_{2}, x^{t+1}_{3}, \dots, x^{2t}_{2}, x^{2t}_{3}$. Thus, $A_{k,t}$ is a mop of order $n = 2(k + 4)t$ and has exactly $n_2 = 3t$ vertices of degree~$2$. We have $n_2 = 3n/(2(k+4))$, satisfying $n/(k+4) < n_2 < 2n/(k+4)$. Clearly, if $S$ is a $K_{1, k + 1}$-isolating set of $A_{k,t}$, then $|S \cap V(F^{i}_{k + 4})| \geq 1$ for each $i \in [2t]$. Thus, $\ii_k(A_{k,t}) \ge 2t$. Since $\{x^{i}_{0} \colon i \in [2t]\}$ is a $K_{1, k + 1}$-isolating set of $A_{k,t}$, we obtain $\ii_k(A_{k,t}) = 2t = n/(k+4)$. 

When $n_2 < n/(k+4)$, we have $\lfloor (n+n_2)/(k+5) \rfloor \le \min\{\lfloor n/(k+4) \rfloor, \lfloor (n-n_2)/(k+2) \rfloor\}$. To see the sharpness of the bound $\lfloor (n+n_2)/(k+5) \rfloor$, let $F^{1}_{k+6},\dots,F^{t}_{k+6}$ be $t$ copies of $F_{k+6}$ and let $F^{1}_{k+4},\dots,F^{t}_{k+4}$ be $t$ copies of $F_{k+4}$, where $F^{1}_{k+6},\dots,F^{t}_{k+6}, F^{1}_{k+4},\dots,F^{t}_{k+4}$ are vertex-disjoint and $(k+4)/2 \le t < k+5$. For $i \in [t]$, let $x^{i}_{0}, x^{i}_{1}, \dots, x^{i}_{k+5}$ be the vertices of $F^i_{k + 6}$ with $x^{i}_{0}$ being the center and $F^i_{k + 6} - \{x^{i}_{0}\} = x^{i}_{1}x^{i}_2 \dots x^{i}_{k+5}$, let $y^{i}_{0}, y^{i}_{1}, \dots, y^{i}_{k+3}$ be the vertices of $F^i_{k + 4}$ with $y^{i}_{0}$ being the center and $F^i_{k + 4} - \{y^{i}_{0}\} = y^{i}_{1}y^{i}_2 \dots y^{i}_{k+3}$, and let $T^{i}_{2k+10}$ be the mop obtained by adding the edges $x^{i}_{k+4}y^{i}_{1}, x^{i}_{k+5}y^{i}_{2}, x^{i}_{k+4}y^{i}_2$ to the union of $F^i_{k + 6}$ and $F^i_{k + 4}$. We extend the union of $T^{1}_{2k+10}, \dots, T^{t}_{2k+10}$ to a mop $H_{k,t}$ by adding edges on the $2t$ vertices $x^{1}_{1}, x^{1}_{2}, \dots, x^{t}_{1}, x^{t}_{2}$. The graph $H_{k,t}$ is illustrated in Figure 1. It is a mop of order $n = (2k + 10)t$ and has exactly $n_2 = t < 2t < (2k+10)t/(k+4) = n/(k+4)$ vertices of degree~$2$. Clearly, if $S$ is a $K_{1, k + 1}$-isolating set of $H_{k,t}$, then, for each $i \in [t]$, $|S \cap V(F^{i}_{k + 6})| \geq 1$ and $|S \cap V(F^{i}_{k + 4})| \geq 1$. Thus, $\ii_k(H_{k,t}) \ge 2t$. Since $\{x_0^i: i \in [t]\} \cup \{y_0^i: i \in [t]\}$ is a $K_{1, k + 1}$-isolating set of $H_{k,t}$, we obtain $\ii_k(H_{k,t}) = 2t$. Since $t < k+5$, $\ii_k(H_{k,t}) = \lfloor 2t + t/(k+5) \rfloor = \lfloor (n + n_2)/(k + 5) \rfloor$. In addition, since $t \ge (k+4)/2$, $\lfloor n/(k+4) \rfloor = \lfloor t(2k+10)/(k+4) \rfloor = \lfloor 2t + 2t/(k+4) \rfloor \ge 2t + 1 > \lfloor (n + n_2)/(k + 5) \rfloor$ and $\lfloor (n-n_2)/(k+2) \rfloor = \lfloor ((2k+10)t-t)/(k+2) \rfloor = \lfloor 2t + 5t/(k+2) \rfloor \ge \lfloor 2t + 5(k+4)/(2(k+2)) \rfloor \ge 2t + 2 > \lfloor (n + n_2)/(k + 5) \rfloor$.
\\

\setlength{\unitlength}{0.8cm}
\begin{center}
\begin{picture}(14,10)

\put(0, 5){\circle*{0.15}}
\put(1, 2.5){\circle*{0.15}}
\put(1, 4.5){\circle*{0.15}}
\put(1.65, 5.8){\circle*{0.15}}
\put(1, 6.5){\circle*{0.15}}
\put(1, 7.5){\circle*{0.15}}
\put(0, 2.5){\circle*{0.15}}
\put(1.5, 3.5){\circle*{0.15}}
\put(1, 4.5){\line(1, 2){0.65}}
\put(1, 6.5){\line(1, -1){0.65}}
\put(1, 4.5){\line(1, -2){0.47}}

\put(1.1, 2.7){\small{$\cdot$}}
\put(1.2, 2.9){\small{$\cdot$}}
\put(1.3, 3.1){\small{$\cdot$}}

\put(0, 2.5){\line(0, 1){2.5}}
\put(0, 2.5){\line(1, 0){13}}
\put(1, 6.5){\line(0, 1){1}}
\put(0, 5){\line(2, 1){1.6}}
\put(0, 5){\line(2, 3){1}}
\put(0, 5){\line(2, 5){1}}
\put(0, 5){\line(2, -1){1}}
\put(0, 5){\line(1, -1){1.5}}
\put(0, 5){\line(2, -5){1}}
\put(2, 6.5){\circle*{0.15}}
\put(2, 7.5){\circle*{0.15}}
\put(2, 8){\circle*{0.15}}

\put(2, 9.5){\circle*{0.15}}
\put(3, 8){\circle*{0.15}}
\put(1, 9){\circle*{0.15}}
\put(2, 6.5){\line(-1, 0){1}}
\put(2, 7.5){\line(-1, 0){1}}
\put(2, 7.5){\line(-1, -1){1}}
\put(2, 6.5){\line(0, 1){1.5}}
\put(2, 9.5){\line(-2, -1){1}}

\put(1.2, 8.2){\small{$\ddots$}}

\put(3, 8){\line(-2, -3){1}}
\put(3, 8){\line(-2, -1){1}}
\put(3, 8){\line(-1, 0){1}}
\put(3, 8){\line(-2, 1){2}}
\put(3, 8){\line(-2, 3){1}}
\put(-0.5, 4.9){\small{$x_0^1$}}
\put(3.2, 7.7){\small{$y_0^1$}}

\put(4, 5){\circle*{0.15}}
\put(5, 2.5){\circle*{0.15}}
\put(5, 4.5){\circle*{0.15}}
\put(5.65, 5.8){\circle*{0.15}}
\put(5, 6.5){\circle*{0.15}}
\put(5, 7.5){\circle*{0.15}}
\put(4, 2.5){\circle*{0.15}}
\put(5.5, 3.5){\circle*{0.15}}
\put(5, 4.5){\line(1, 2){0.65}}
\put(5, 6.5){\line(1, -1){0.65}}
\put(5, 4.5){\line(1, -2){0.47}}
\put(5.1, 2.7){\small{$\cdot$}}
\put(5.2, 2.9){\small{$\cdot$}}
\put(5.3, 3.1){\small{$\cdot$}}
\put(4, 2.5){\line(0, 1){2.5}}
\put(4, 2.5){\line(1, 0){1}}
\put(5, 6.5){\line(0, 1){1}}
\put(4, 5){\line(2, 1){1.6}}
\put(4, 5){\line(2, 3){1}}
\put(4, 5){\line(2, 5){1}}
\put(4, 5){\line(2, -1){1}}
\put(4, 5){\line(1, -1){1.5}}
\put(4, 5){\line(2, -5){1}}
\put(6, 6.5){\circle*{0.15}}
\put(6, 7.5){\circle*{0.15}}
\put(6, 8){\circle*{0.15}}

\put(6, 9.5){\circle*{0.15}}
\put(7, 8){\circle*{0.15}}
\put(5, 9){\circle*{0.15}}
\put(6, 6.5){\line(-1, 0){1}}
\put(6, 7.5){\line(-1, 0){1}}
\put(6, 7.5){\line(-1, -1){1}}
\put(6, 6.5){\line(0, 1){1.5}}
\put(6, 9.5){\line(-2, -1){1}}
\put(5.2, 8.2){\small{$\ddots$}}
\put(7, 8){\line(-2, -3){1}}
\put(7, 8){\line(-2, -1){1}}
\put(7, 8){\line(-1, 0){1}}
\put(7, 8){\line(-2, 1){2}}
\put(7, 8){\line(-2, 3){1}}
\put(3.5, 4.9){\small{$x_0^2$}}
\put(7.2, 7.7){\small{$y_0^2$}}

\put(12, 5){\circle*{0.15}}
\put(13, 2.5){\circle*{0.15}}
\put(13, 4.5){\circle*{0.15}}
\put(13.65, 5.8){\circle*{0.15}}
\put(13, 6.5){\circle*{0.15}}
\put(13, 7.5){\circle*{0.15}}
\put(12, 2.5){\circle*{0.15}}
\put(13.5, 3.5){\circle*{0.15}}
\put(13, 4.5){\line(1, 2){0.65}}
\put(13, 6.5){\line(1, -1){0.65}}
\put(13, 4.5){\line(1, -2){0.47}}
\put(13.1, 2.7){\small{$\cdot$}}
\put(13.2, 2.9){\small{$\cdot$}}
\put(13.3, 3.1){\small{$\cdot$}}
\put(12, 2.5){\line(0, 1){2.5}}
\put(12, 2.5){\line(1, 0){1}}
\put(13, 6.5){\line(0, 1){1}}
\put(12, 5){\line(2, 1){1.6}}
\put(12, 5){\line(2, 3){1}}
\put(12, 5){\line(2, 5){1}}
\put(12, 5){\line(2, -1){1}}
\put(12, 5){\line(1, -1){1.5}}
\put(12, 5){\line(2, -5){1}}
\put(14, 6.5){\circle*{0.15}}
\put(14, 7.5){\circle*{0.15}}
\put(14, 8){\circle*{0.15}}

\put(14, 9.5){\circle*{0.15}}
\put(15, 8){\circle*{0.15}}
\put(13, 9){\circle*{0.15}}
\put(14, 6.5){\line(-1, 0){1}}
\put(14, 7.5){\line(-1, 0){1}}
\put(14, 7.5){\line(-1, -1){1}}
\put(14, 6.5){\line(0, 1){1.5}}
\put(14, 9.5){\line(-2, -1){1}}
\put(13.2, 8.2){\small{$\ddots$}}
\put(15, 8){\line(-2, -3){1}}
\put(15, 8){\line(-2, -1){1}}
\put(15, 8){\line(-1, 0){1}}
\put(15, 8){\line(-2, 1){2}}
\put(15, 8){\line(-2, 3){1}}
\put(11.5, 4.9){\small{$x_0^t$}}
\put(15.2, 7.7){\small{$y_0^t$}}

\put(9, 4.9){\small{$...$}}
\put(9, 7.7){\small{$...$}}
\put(9, 2){\small{$...$}}

\qbezier(0, 2.5)(6.5, -3.5)(13, 2.5)
\qbezier(1, 2.5)(7, -2.5)(13, 2.5)
\qbezier(4, 2.5)(8.5, -1)(13, 2.5)
\qbezier(4, 2.5)(8, -0.5)(12, 2.5)
\qbezier(5, 2.5)(8.5, -0.2)(12, 2.5)

\end{picture}
\end{center}
\begin{center} \vskip 25 pt \footnotesize{\textbf{Figure 1 :} $H_{k,t}$}
\end{center}
\vskip 30 pt

We now show that the actual bound $(n + n_2)/(k + 5)$ in Theorem~\ref{t:main2} is attained for $n_2 = n/(k+4)$. For an integer $t \ge 1$, let $F^1_{k + 4}, \dots, F^{t}_{k + 4}$ be $t$ vertex-disjoint fans of order $k+4$. For $i \in [t]$, let $x^{i}_{0}, x^{i}_{1}, \dots, x^{i}_{k+3}$ be the vertices of $F^i_{k + 4}$ with $x^{i}_{0}$ being the center and with $x^{i}_{1}$ and $x^{i}_{k+3}$ being the vertices of degree $2$. We extend the union of $F^1_{k + 4}, \dots, F^{t}_{k + 4}$ to a mop $T_{k,t}$ by adding edges on the $2t$ vertices $x^{1}_{1}, x^{1}_{2}, \dots, x^{t}_{1}, x^{t}_{2}$. Thus, the order $n$ of $T_{k,t}$ is $(k + 4)t$, and $x^{1}_{k+3}, \dots, x^{t}_{k+3}$ are the vertices of $T_{k,t}$ of degree~$2$. Clearly, if $S$ is a $K_{1, k + 1}$-isolating set of $T_{k,t}$, then $|S \cap V(F^{i}_{k + 4})| \geq 1$ for each $i \in [t]$. Thus, $\ii_k(T_{k,t}) \ge t$. Since $\{x^{i}_{0} \colon i \in [t]\}$ is a $K_{1, k + 1}$-isolating set of $T_{k,t}$, $\ii_k(T_{k,t}) = t$. Since $n = (k + 4)t$ and $n_2 = t$, it follows that $\ii_k(T_{k,t}) = n/(k+4) = t = (n+n_2)/(k+5)$.

When $n_2 > 2n/(k+4)$, we have $(n-n_2)/(k+2) < \min\{ n/(k+4), (n+n_2)/(k+5)\}$. We demonstrate that the bound $(n-n_2)/(k+2)$ is sharp. Let $x$ be the center of $F_{k+2}$ and let $y_1, \dots, y_{k+1}$ be $k+1$ distinct isolated vertices. Let $uv$ be a Hamiltonian edge of $F_{k+2}$ with $d_{F_{k+2}}(u) = d_{F_{k+2}}(v) = 3$, and let $z_1z_1', \dots, z_{k+1}z_{k+1}'$ be the remaining Hamiltonian edges of $F_{k+2}$. Let $R_{2k+3}$ be the graph with vertex set $V(F_{k+2}) \cup \{y_1, \dots, y_{k+1}\}$ and edge set $E(F_{k+2}) \cup \{y_1z_1, y_1z_1', \dots, y_{k+1}z_{k+1}, y_{k+1}z_{k+1}'\}$. Clearly, $R_{2k+3}$ is a $(2k+3)$-vertex mop with exactly $k+1$ vertices of degree $2$. Let $R^{1}_{2k+3}, \dots, R^{t}_{2k+3}$ be $t$ vertex-disjoint copies of $R_{2k+3}$. For $i \in [t]$, let $x_i$, $u_i$, and $v_i$ be the vertices of $R^{i}_{2k+3}$ corresponding to the vertices $x$, $u$, and $v$ of $R_{2k+3}$, respectively. We extend the union of $R^{1}_{2k+3}, \dots, R^{t}_{2k+3}$ to a mop $S_{k,t}$ by adding edges on the $2t$ vertices $u_1, v_1, \dots, u_t, v_t$ (see Figure 2). Thus, the order $n$ of $S_{k,t}$ is $(2k+3)t$, and the number $n_2$ of vertices of $S_{k,t}$ of degree~$2$ is $(k+1)t$. We have $n_2 = (k+1)n/(2k+3) > 2n/(k+4)$ as $k \geq 2$. Clearly, if $D$ is a $K_{1, k + 1}$-isolating set of $S_{k,t}$, then $|D \cap V(R^{i}_{2k + 3})| \geq 1$ for each $i \in [t]$. Thus, $\ii_k(S_{k,t}) \ge t$. Since $\{x_i \colon i \in [t]\}$ is a $K_{1, k + 1}$-isolating set of $S_{k,t}$, $\ii_k(S_{k,t}) = t = (n - n_2)/(k + 2)$. 

\setlength{\unitlength}{0.9cm}
\begin{center}
\begin{picture}(12,8)

\put(1, 6.7){\small{$x_{1}$}}
\put(1, 6.5){\circle*{0.15}}
\put(0, 5.5){\circle*{0.15}}
\put(0, 4.5){\circle*{0.15}}
\put(0, 3.5){\circle*{0.15}}
\put(0, 2.5){\circle*{0.15}}
\put(2, 5.5){\circle*{0.15}}
\put(2, 4.5){\circle*{0.15}}
\put(2, 3.5){\circle*{0.15}}
\put(2, 2.5){\circle*{0.15}}
\put(-0.2, 2.2){\small{$u_{1}$}}
\put(1.8, 2.2){\small{$v_{1}$}}
\put(-1, 6.5){\circle*{0.15}}
\put(-1, 5.5){\circle*{0.15}}
\put(-1, 4.5){\circle*{0.15}}
\put(-1, 3.5){\circle*{0.15}}
\put(3, 6.5){\circle*{0.15}}
\put(3, 5.5){\circle*{0.15}}
\put(3, 4.5){\circle*{0.15}}
\put(3, 3.5){\circle*{0.15}}
\put(1, 6.5){\line(1, -1){1}}
\put(1, 6.5){\line(1, -2){1}}
\put(1, 6.5){\line(1, -3){1}}
\put(1, 6.5){\line(1, -4){1}}
\put(1, 6.5){\line(-1, -1){1}}
\put(1, 6.5){\line(-1, -2){1}}
\put(1, 6.5){\line(-1, -3){1}}
\put(1, 6.5){\line(-1, -4){1}}

\put(0, 2.5){\line(0, 1){3}}
\put(2, 2.5){\line(0, 1){3}}
\put(-1, 3.5){\line(1, 0){1}}
\put(-1, 3.5){\line(1, -1){1}}
\put(-1, 4.5){\line(1, 0){1}}
\put(-1, 4.5){\line(1, -1){1}}
\put(-1, 5.5){\line(1, 0){1}}
\put(-1, 5.5){\line(1, -1){1}}
\put(-1, 6.5){\line(1, 0){2}}
\put(-1, 6.5){\line(1, -1){1}}

\put(3, 3.5){\line(-1, 0){1}}
\put(3, 3.5){\line(-1, -1){1}}
\put(3, 4.5){\line(-1, 0){1}}
\put(3, 4.5){\line(-1, -1){1}}
\put(3, 5.5){\line(-1, 0){1}}
\put(3, 5.5){\line(-1, -1){1}}
\put(3, 6.5){\line(-1, 0){2}}
\put(3, 6.5){\line(-1, -1){1}}
\put(0, 2.5){\line(1, 0){12}}

\put(6, 6.7){\small{$x_{2}$}}
\put(6, 6.5){\circle*{0.15}}
\put(5, 5.5){\circle*{0.15}}
\put(5, 4.5){\circle*{0.15}}
\put(5, 3.5){\circle*{0.15}}
\put(5, 2.5){\circle*{0.15}}
\put(7, 5.5){\circle*{0.15}}
\put(7, 4.5){\circle*{0.15}}
\put(7, 3.5){\circle*{0.15}}
\put(7, 2.5){\circle*{0.15}}
\put(4.8, 2.2){\small{$u_{2}$}}
\put(6.8, 2.2){\small{$v_{2}$}}
\put(4, 6.5){\circle*{0.15}}
\put(4, 5.5){\circle*{0.15}}
\put(4, 4.5){\circle*{0.15}}
\put(4, 3.5){\circle*{0.15}}
\put(8, 6.5){\circle*{0.15}}
\put(8, 5.5){\circle*{0.15}}
\put(8, 4.5){\circle*{0.15}}
\put(8, 3.5){\circle*{0.15}}
\put(6, 6.5){\line(1, -1){1}}
\put(6, 6.5){\line(1, -2){1}}
\put(6, 6.5){\line(1, -3){1}}
\put(6, 6.5){\line(1, -4){1}}
\put(6, 6.5){\line(-1, -1){1}}
\put(6, 6.5){\line(-1, -2){1}}
\put(6, 6.5){\line(-1, -3){1}}
\put(6, 6.5){\line(-1, -4){1}}

\put(5, 2.5){\line(0, 1){3}}
\put(7, 2.5){\line(0, 1){3}}
\put(4, 3.5){\line(1, 0){1}}
\put(4, 3.5){\line(1, -1){1}}
\put(4, 4.5){\line(1, 0){1}}
\put(4, 4.5){\line(1, -1){1}}
\put(4, 5.5){\line(1, 0){1}}
\put(4, 5.5){\line(1, -1){1}}
\put(4, 6.5){\line(1, 0){2}}
\put(4, 6.5){\line(1, -1){1}}

\put(8, 3.5){\line(-1, 0){1}}
\put(8, 3.5){\line(-1, -1){1}}
\put(8, 4.5){\line(-1, 0){1}}
\put(8, 4.5){\line(-1, -1){1}}
\put(8, 5.5){\line(-1, 0){1}}
\put(8, 5.5){\line(-1, -1){1}}
\put(8, 6.5){\line(-1, 0){2}}
\put(8, 6.5){\line(-1, -1){1}}
\put(5, 2.5){\line(1, 0){2}}

\put(11, 6.7){\small{$x_{3}$}}
\put(11, 6.5){\circle*{0.15}}
\put(10, 5.5){\circle*{0.15}}
\put(10, 4.5){\circle*{0.15}}
\put(10, 3.5){\circle*{0.15}}
\put(10, 2.5){\circle*{0.15}}
\put(12, 5.5){\circle*{0.15}}
\put(12, 4.5){\circle*{0.15}}
\put(12, 3.5){\circle*{0.15}}
\put(12, 2.5){\circle*{0.15}}
\put(10, 2.2){\small{$u_{3}$}}
\put(12.2, 2.2){\small{$v_{3}$}}
\put(9, 6.5){\circle*{0.15}}
\put(9, 5.5){\circle*{0.15}}
\put(9, 4.5){\circle*{0.15}}
\put(9, 3.5){\circle*{0.15}}
\put(13, 6.5){\circle*{0.15}}
\put(13, 5.5){\circle*{0.15}}
\put(13, 4.5){\circle*{0.15}}
\put(13, 3.5){\circle*{0.15}}
\put(11, 6.5){\line(1, -1){1}}
\put(11, 6.5){\line(1, -2){1}}
\put(11, 6.5){\line(1, -3){1}}
\put(11, 6.5){\line(1, -4){1}}
\put(11, 6.5){\line(-1, -1){1}}
\put(11, 6.5){\line(-1, -2){1}}
\put(11, 6.5){\line(-1, -3){1}}
\put(11, 6.5){\line(-1, -4){1}}

\put(10, 2.5){\line(0, 1){3}}
\put(12, 2.5){\line(0, 1){3}}
\put(9, 3.5){\line(1, 0){1}}
\put(9, 3.5){\line(1, -1){1}}
\put(9, 4.5){\line(1, 0){1}}
\put(9, 4.5){\line(1, -1){1}}
\put(9, 5.5){\line(1, 0){1}}
\put(9, 5.5){\line(1, -1){1}}
\put(9, 6.5){\line(1, 0){2}}
\put(9, 6.5){\line(1, -1){1}}

\put(13, 3.5){\line(-1, 0){1}}
\put(13, 3.5){\line(-1, -1){1}}
\put(13, 4.5){\line(-1, 0){1}}
\put(13, 4.5){\line(-1, -1){1}}
\put(13, 5.5){\line(-1, 0){1}}
\put(13, 5.5){\line(-1, -1){1}}
\put(13, 6.5){\line(-1, 0){2}}
\put(13, 6.5){\line(-1, -1){1}}
\put(10, 2.5){\line(1, 0){2}}

\qbezier(0, 2.5)(6, -1.5)(12, 2.5)
\qbezier(0, 2.5)(5, -1)(10, 2.5)
\qbezier(2, 2.5)(6, -0.5)(10, 2.5)
\qbezier(5, 2.5)(7.5, 0.5)(10, 2.5)

\end{picture}
\end{center}
\begin{center} \footnotesize{\textbf{Figure 2 :} $S_{k,t}$ with $k = 7$ and $t = 3$}
\end{center}
\vskip 30 pt

\section{Proofs of Theorems~\ref{t:main1}--\ref{t:MAIN}} \label{S:main1}

In this section, we prove Theorems~\ref{t:main1}--\ref{t:MAIN}. We apply results of O'Rourke~\cite{Ro83} in computational geometry that were used in a new proof by Lema\'{n}ska, Zuazua and Zylinski~\cite{LeZuZy17} of an upper bound by Dorfling, Hattingh and Jonck \cite{DoHaJo16} on the size of a total dominating set~(a set $S$ of vertices such that each vertex of the graph is adjacent to a vertex in $S$) of a mop. Before stating these results, we make a related straightforward observation that we will also use.

Given three mops $G$, $G_1$ and $G_2$, we say that a diagonal $d$ of $G$ \emph{partitions $G$ into $G_1$ and $G_2$} if $G$ is the union of $G_1$ and $G_2$, $V(G_1) \cap V(G_2) = d$, and $E(G_1) \cap E(G_2) = \{d\}$. The following is a straightforward well-known fact.

\parskip=0pt

\begin{lem}\label{diagonallemma} If $d$ is a diagonal of a mop $G$, then $d$ partitions $G$ into two mops $G_1$ and $G_2$.
\end{lem}


\begin{lem}
\label{mop:6-10}
If $r \ge 0$ and $G$ is a mop of order $n \geq 2r+4$, then $G$ has a diagonal $d$ that partitions it into two mops $G_1$ and $G_2$ such that $G_1$ has exactly $\ell$ Hamiltonian edges of $G$ for some $\ell \in [r+2, 2r+2]$.
\end{lem}
\proof
In this proof, all subscripts are taken modulo $n$. Let $x_0x_1 \dots x_{n-1}x_0$ be the Hamiltonian cycle of $G$. Let $p$ be the smallest integer such that $p \geq r+2$ and $x_{i}x_{i+p} \in E(G)$ for some $i \in \{0, 1, \dots, n-1\}$ ($p$ exists as $x_0x_{n-1} \in E(G)$). It suffices to show that $p \le 2r + 2$. For some $q \in [p-1]$, $G$ has a triangular face containing $x_i$, $x_{i+q}$, and $x_{i+p}$, so $x_ix_{i+q},x_{i+q}x_{i+p} \in E(G)$. By the choice of $p$, we have $q \le r + 1$ and $p - q = (i+p) - (i+q) \leq r + 1$. We have $p = p - q + q \leq 2r + 2$.~\hfill{\qed}

\parskip=4mm
The case $r = 2$ of Lemma~\ref{mop:6-10} was proved by Chv\'{a}tal in \cite{Ch75} and is restated in \cite[Lemma~1.1]{Ro87}. 
The case $r = 3$ of Lemma~\ref{mop:6-10} was proved by O'Rourke \cite{Ro83}.

\parskip=4pt
For a graph $G$ and an edge $uv$ of $G$, the \emph{edge contraction} of $G$ along $uv$ is the graph obtained from $G$ by deleting $u$ and $v$ (and all incident edges), adding a new vertex $x$, and making $x$ adjacent to the vertices in $N_G[\{u,v\}] \setminus \{u,v\}$ only. Recall that every mop can be embedded in a plane so that the exterior face contains all vertices. By looking at polygon corners as vertices, we have that a mop is a triangulation of a simple polygon, meaning that its boundary is the polygon and its interior faces are triangles.

\parskip=0pt
\begin{lem}{\rm (\cite{Ro83})}
\label{contraction}
If $G$ is a triangulation of order at least $4$ of a simple polygon $P$, $e$ is a Hamiltonian edge of $G$, and $G'$ is the edge contraction of $G$ along $e$, then $G'$ is a triangulation of some simple polygon $P'$.
\end{lem}

If $G$ is a graph and $I \subseteq V(G)$ such that $uv \notin E(G)$ for every $u, v \in I$, then $I$ is called an \emph{independent set of $G$}.




\parskip=0pt

\begin{lem}\label{maxofn2}
If $G$ is a mop of order $n \geq 4$, then the set of vertices of $G$ of degree $2$ is an independent set of $G$ of size at most $\frac{n}{2}$.
\end{lem}
Lemma~\ref{maxofn2} is given in \cite{BK}, as are parts (a)--(d) of the next lemma. 

\begin{lem}\label{basicfacts} If $G$ is a mop of order $n \geq 3$, then the following assertions hold:\\
(a) Each vertex of $G$ is of degree at least $2$. \\
(b) $G$ has at least $2$ vertices of degree $2$. \\
(c) If $n \geq 4$, then $G - v$ is a mop for each vertex $v$ of $G$ of degree $2$. \\
(d) A graph $H$ is a mop if $G = H-w$ for some $w \in V(H)$ such that $d_H(w) = 2$ and $N_H(w)$ is a Hamiltonian edge of $G$.\\
(e) Each vertex of $G$ has at most two neighbours of degree $2$.
\end{lem}

\proof For (a)--(d), see \cite{BK}. We prove (e). Let $V_2$ be the set of vertices of $G$ of degree $2$. Let $v \in V(G)$. We may label the vertices $x_1, x_2, \dots, x_n$ so that $x_1x_2 \dots x_nx_1$ is the Hamiltonian cycle of $G$ and $x_n = v$. If $i \in [2,n-2]$ such that $x_i \in N_G(v)$, then $x_{i-1}, x_{i+1}, v \in N_G(x_i)$, so $x_i \notin V_2$. Thus, $N_G(v) \cap V_2 \subseteq \{x_1, x_{n-1}\}$.~\hfill{\qed}
\parskip=4mm

Lemmas~\ref{maxofn2} and \ref{basicfacts}(b) tell us that the number $n_2$ of vertices of degree $2$ (of a mop) satisfies
\begin{equation} 2 \leq n_2 \leq \frac{n}{2}. \label{n_2bounds}
\end{equation}
We mention that both bounds are sharp \cite{BK}.

\parskip=4pt

The next lemma settles Theorem~\ref{t:main1} for $n \leq 2k+7$, and hence allows us to use Lemma~\ref{mop:6-10} in the proof of Theorem~\ref{t:main1}. 

We say that a vertex $x$ of a mop $G$ is a \emph{diagonal $K_{1,k+1}$-isolating vertex of $G$} if $\{x\}$ is a $K_{1,k+1}$-isolating set of $G$ and $x$ is one of the two vertices of a diagonal of $G$.

\begin{lem}\label{mop:iotaofsmallmop}
If $k \ge 0$ and $G$ is a mop of order $n \le 2k + 7$, then 
$\ii_{k}(G) \le 1$.
\end{lem}
\proof Suppose $n = 2k+7$. Let $x_1x_2 \dots x_{2k+7}x_1$ be the unique Hamiltonian cycle of $G$ and hence the boundary of the exterior face of $G$. Let $r = k+1$. Thus, $n > 2r+4$. By Lemma~\ref{mop:6-10}, $G$ has a diagonal $d$ that partitions it into two mops $G_1$ and $G_2$ such that $G_1$ has exactly $\ell$ Hamiltonian edges of $G$ for some $\ell \in [r+2, 2r+2] = [k+3,2k+4]$. We may assume that $d$ is the edge $x_1x_{\ell+1}$ and that $V(G_1) = \{x_1, x_2, \dots, x_{\ell+1}\}$. Thus, $V(G_2) = \{x_1, x_{\ell + 1},x_{\ell + 2},\dots,x_{2k+7}\}$.
\parskip=4pt

Suppose $\ell = k+3$. Then, $|V(G_1)| = k+4$, $|V(G_2)| = k+5$, $x_{k+4}, x_n \in N_{G_2}(x_1)$, and $x_1, x_{k+5} \in N_{G_2}(x_{k+4})$. 
Since $x_1$ and $x_{k+4}$ are adjacent in $G_2$, Lemma~\ref{maxofn2} tells us that their degrees in $G_2$ cannot be both $2$. Thus, $|N_{G_2}[x_i]| \ge 4$ for some $i \in \{1,k+4\}$. Clearly, $|N_{G_1}[x_i]| \ge 3$. Thus, $|V(G_1) \setminus N_{G}[x_i]| \le k+1$ and $|V(G_2) \setminus N_{G}[x_i]| \le k+1$, and hence $G - N_{G}[x_i]$ contains no copy of $K_{1,k+1}$. Therefore, $x_i$ is a diagonal $K_{1,k+1}$-isolating vertex of $G$.

We now proceed by induction on $\ell$. 
Thus, we consider $\ell \ge k + 4$ and assume that, if $G$ has a diagonal that partitions it into two mops $H_1$ and $H_2$ such that $H_1$ has exactly $\ell^*$ Hamiltonian edges of $G$ for some $\ell^* \in [k+3, \ell - 1]$, then $G$ has a diagonal $K_{1,k+1}$-isolating vertex. Let $(x_1,x_i,x_{\ell + 1})$ be the triangular face of $G_1$ containing the Hamiltonian edge $x_1x_{\ell + 1}$ of $G_1$. Thus, $2 \le i \le \ell$. Let $\ell' = \ell+1-i$ and $\ell'' = i-1$. Suppose $k+4 \le i \le \ell$. By Lemma~\ref{diagonallemma}, $x_1x_i$ partitions $G$ into two mops $G'_1$ and $G'_2$ such that $G'_1$ contains the $\ell''$ Hamiltonian edges $x_1x_2,x_2x_3,\dots,x_{i-1}x_{i}$ of $G$. Since $k+4 \le i \le \ell$, $k+3 \le \ell'' \le \ell-1$. By the induction hypothesis, $G$ has a diagonal $K_{1,k+1}$-isolating vertex. Now suppose $2 \le i \le k+3$. By Lemma~\ref{diagonallemma}, $x_ix_{\ell+1}$ partitions $G$ into two mops $G'_1$ and $G'_2$ such that $G'_1$ contains the $\ell'$ Hamiltonian edges $x_ix_{i+1},x_{i+1}x_{i+2},\dots,x_{\ell}x_{\ell+1}$ of $G$. If $2 \le i \le \ell-k-2$, then $k+3 \le \ell' \le \ell-1$, and hence, by the induction hypothesis, $G$ has a diagonal $K_{1,k+1}$-isolating vertex. Suppose $\ell-k-1 \le i \le k+3$. Since $x_1x_2,x_1x_{\ell+1} \in E(G_1)$, it follows that $|\{x_1,x_2,\dots,x_{i-1}\} \setminus N_{G_1}[x_{1}]| \le i-3 \le (k+3)-3 = k$ and $|\{x_{i+1},x_{i+2},\dots,x_{\ell+1}\} \setminus N_{G_1}[x_{1}]| \le \ell-i \le \ell-(\ell-k-1) = k+1$. Thus, $G_1 - N_{G_1}[x_{1}]$ contains no copy of $K_{1,k+1}$. Now $|V(G_2)| = 2k+8- \ell \le k+4$ as $\ell \ge k+4$. Since $x_1, x_{\ell + 1}, x_n \in N_{G_2}[x_1]$, $|V(G_2) \setminus N_{G_2}[x_1]| \le k+1$. Thus, $G_2 - N_{G_2}[x_{1}]$ contains no copy of $K_{1,k+1}$. Therefore, $x_{1}$ is a diagonal $K_{1,k+1}$-isolating vertex of $G$.

Now suppose $n = 2k+6$. Let $uv$ be a Hamiltonian edge of $G$. By Lemma~\ref{basicfacts}(d), we can obtain a mop $H$ from $G$ by inserting a vertex $w$ in the exterior face of $G$ and adding the edges $wu$ and $wv$. Since $H$ is a mop of order $2k+7$, $H$ has a diagonal $K_{1,k+1}$-isolating vertex $x$. Since $w$ is not a vertex of a diagonal of $H$, $x \neq w$. Thus, $\{x\}$ is a $K_{1,k+1}$-isolating set of $G$.

For $i \leq 2k+5$, we obtain the result for $n = i$ from the result for $n = i+1$ in the same way we obtained the result for $n = 2k+6$ from the result for $n = 2k+7$.~\hfill{\qed}

\parskip=4mm

We now prove Theorems~\ref{t:main1}--\ref{t:MAIN}. Recall that the bounds in Theorems~\ref{t:main1}--\ref{t:main3} are sharp by the constructions in Section~\ref{extremalsection}, so we now prove the bounds.  

\parskip=4mm


\noindent \textbf{Proof of Theorem~\ref{t:main1}.} 
If $k+4 \le n \le 2k+7$, then $\ii_k(G) \le 1 \leq n/(k+4)$ by Lemma~\ref{mop:iotaofsmallmop}. We now assume that $n \ge 2k + 8$ and proceed by induction on $n$. Let $x_1x_2 \ldots x_nx_1$ be the unique Hamiltonian cycle $C$ of $G$ and hence the boundary of the exterior face of $G$. By Lemma~\ref{mop:6-10} with $r = k+2$, $G$ has a diagonal $d$ that partitions it into two mops $G_1$ and $G_2$ such that $G_1$ has exactly $\ell$ Hamiltonian edges of $G$ for some $\ell \in [k+4, 2k+6]$. We may assume that $d = x_1x_{\ell + 1}$ and $V(G_1) =  \{x_1,x_2,\ldots,x_{\ell + 1}\}$. Note that $x_1x_2,x_2x_3,\dots,x_{\ell}x_{\ell + 1}$ are the $\ell$ Hamiltonian edges of $G$ that belong to $G_1$. Let $(x_1,x_j,x_{\ell + 1})$ be the triangular face of $G_1$ that contains the edge $x_1x_{\ell + 1}$. Then, $2 \le j \le \ell$. 

\parskip=4pt
Suppose $\ell = k + 4$. Let $G'$ be the graph obtained from $G$ by deleting the vertices $x_2,x_3,\dots,x_{k + 4}$ and contracting the edge $x_1x_{k + 5}$ to form a new vertex $z$ (see Figure 3). Thus, $G'$ is obtained from $G_2$ by contracting the edge $x_1x_{k + 5}$. By Lemma  \ref{contraction}, $G'$ is a mop. Let $n' = |V(G')|$. Thus, $n' = n - (k + 4) \ge k + 4$. By the induction hypothesis, $\ii_k(G') \le n'/(k + 4) = n/(k + 4)-1$. Let $S'$ be a smallest $K_{1, k + 1}$-isolating set of $G'$. Then, $|S'| = \ii_k(G') \le n/(k + 4)-1$. Since $|V(G_1) \setminus N_{G_1}[\{x_1,x_{k+5}\}]| \le k+1$, $\{x_1,x_{k+5}\}$ is a $K_{1,k+1}$-isolating set of $G_1$. We have $N_{G'}(z) \subseteq N_{G}[\{x_1, x_{k+5}\}]$. Thus, if $z \in S'$, then $(S' \setminus \{z\}) \cup \{x_1,x_{k+5}\}$ is a $K_{1,k+1}$-isolating set of $G$, and hence $\ii_k(G) \le (|S'| - 1) + 2 \le n/(k+4)$. Suppose $z \notin S'$. Since $x_1, x_{j-1}, x_j, x_{j+1}, x_{k+5} \in N_{G_1}[x_j]$ and $\{x_{j-1}, x_{j+1}\} \setminus \{x_1, x_{k+5}\} \neq \emptyset$ (as $2 \le j \le \ell = k+4$), we have $|V(G_1) \setminus N_{G_1}[x_j]| \le k+1$, so $\{x_j\}$ is a $K_{1,k+1}$-isolating set of $G_1$. Since $x_1, x_{k + 5} \in N_G[x_j]$, it follows that $S' \cup \{x_j\}$ is a $K_{1, k + 1}$-isolating set of $G$. Therefore, $\ii_k(G) \le |S'| + 1 \le n/(k + 4)$.\medskip

\setlength{\unitlength}{0.6cm}
\begin{center}
\begin{picture}(15, 6)
\put(-1, 0){\circle*{0.2}}
\put(3, 0){\circle*{0.2}}
\put(-2, 2){\circle*{0.2}}
\put(-1, 4){\circle*{0.2}}
\put(3, 4){\circle*{0.2}}
\put(4, 2){\circle*{0.2}}
\put(1, 5){\circle*{0.2}}
\put(-4, 0.2){\small{$x_n$}}
\put(6, 0.2){\small{$x_{k + 6}$}}
\put(-4, -3.5){\small{$x_p$}}
\put(0, -3.5){\small{$x_q$}}
\put(6, -3.5){\small{$x_r$}}
\put(-4, 0){\line(1, 0){10}}

\put(-1, 0){\line(-1, 2){1}}
\put(-1, 4){\line(-1, -2){1}}
\put(-1, 4){\line(2, 1){2}}

\put(4, 2){\line(-1, 2){1}}
\put(4, 2){\line(-1, -2){1}}
\put(-2, 0.2){\small{$x_1$}}
\put(-2.8, 2){\small{$x_2$}}
\put(-1.4, 4.2){\small{$x_3$}}
\put(1, 5.2){\small{$x_4$}}
\put(1.8, 4.4){\small{$\cdot$}}
\put(2, 4.3){\small{$\cdot$}}
\put(2.2, 4.2){\small{$\cdot$}}
\put(3, 4.2){\small{$x_{k + 3}$}}
\put(4.2, 2){\small{$x_{k + 4}$}}
\put(3.4, 0.2){\small{$x_{k + 5}$}}

\put(-4, 0){\circle*{0.2}}
\put(6, 0){\circle*{0.2}}
\put(1, -3){\circle*{0.2}}
\put(-4, -3){\circle*{0.2}}
\put(6, -3){\circle*{0.2}}
\put(6, -3){\line(-1, 1){3}}
\put(-4, -3){\line(1, 1){3}}
\put(1, -3){\line(-2, 3){2}}
\put(1, -3){\line(2, 3){2}}

\put(-1.7, -2){\small{$...$}}
\put(3, -2){\small{$...$}}
\put(10, 0.2){\small{$x_n$}}
\put(18, 0.2){\small{$x_{k + 6}$}}
\put(14, 0.2){\small{$z$}}
\put(10, -3.5){\small{$x_p$}}
\put(18, -3.5){\small{$x_r$}}
\put(14, -3.5){\small{$x_q$}}

\put(14, 0){\circle*{0.2}}
\put(10, 0){\circle*{0.2}}
\put(18, 0){\circle*{0.2}}
\put(14, -3){\circle*{0.2}}
\put(10, -3){\circle*{0.2}}
\put(18, -3){\circle*{0.2}}
\put(18, -3){\line(-4, 3){4}}
\put(10, -3){\line(4, 3){4}}
\put(14, -3){\line(0, 1){3}}
\put(10, 0){\line(1, 0){8}}
\put(12, -2){\small{$...$}}
\put(15, -2){\small{$...$}}
\put(7.5, 0){\Large{$\Rightarrow$}}

\end{picture}
\end{center} 
\vskip 75 pt
\begin{center} \footnotesize{\textbf{\textbf{Figure 3 :}} The edge contraction of $G_2$.}
\end{center}
\vskip 30 pt


We now use induction on $\ell$. Thus, we consider $\ell \ge k + 5$ and assume that, if $G$ has a diagonal that partitions it into two mops $H_1$ and $H_2$ such that $H_1$ has exactly $\ell^*$ Hamiltonian edges of $G$ for some $\ell^* \in [k+4, \ell - 1]$, then $\ii_k(G) \leq n/(k+4)$. Since $\ell \le 2k+6$, $\ell = k + 4 + t$ for some $t \in [k + 2]$.

\parskip=0pt

\begin{claim}
\label{c1}
If $j \notin [t + 2, k + 4]$, then $\ii_k(G) \leq n/(k+4)$.
\end{claim}
\proof Let $\ell' = j-1$ and $\ell'' = \ell + 1 - j$. Suppose $k + 5 \le j \le \ell$. Then, $k+4 \le \ell' \le \ell - 1$. By Lemma~\ref{diagonallemma}, $x_1x_{j}$ partitions $G$ into two mops $G'_1$ and $G'_2$ such that $G'_1$ contains the $\ell'$ Hamiltonian edges $x_1x_2,x_2x_3,\dots,x_{j-1}x_j$ of $G$. By the induction hypothesis, $\ii_k(G) \leq n/(k+4)$. Now suppose $2 \le j \le t + 1$.  Since $\ell = k + 4 + t$, $k + 4 \le \ell'' \le \ell - 1$. By Lemma~\ref{diagonallemma}, $x_jx_{\ell}$ partitions $G$ into two mops $G''_1$ and $G''_2$ such that $G''_1$ contains the $\ell''$ Hamiltonian edges $x_jx_{j + 1},x_{j + 1}x_{j + 2},\dots,x_{\ell}x_{\ell + 1}$ of $G$. By the induction hypothesis, $\ii_k(G) \leq n/(k+4)$.\hfill{\smallqed}

\parskip=3mm

In view of Claim~\ref{c1}, we now assume that $j \in [t + 2, k + 4]$.

\parskip=0pt

\begin{claim}
\label{c2}
$G_1 - N_{G_{1}}[x_j]$ contains no copy of $K_{1, k + 1}$.
\end{claim}
\proof Since $x_1, x_{j-1}, x_{j+1}, x_{\ell + 1} \in N_G(x_{j})$ and $j \in [t + 2, k + 4]$, there are at most $j - 3 \le k + 1$ vertices in $\{x_1,x_2,\dots,x_{j - 1}\}$ which are not adjacent to $x_{j}$, and at most $\ell + 1 - (j + 2) \le (k + 4 + t + 1) - (t + 4) \le k + 1$ vertices in $\{x_{j+1},x_{j+2},\dots,x_{\ell + 1}\}$ which are not adjacent to $x_{j}$. Since no vertex in $\{x_2,\dots,x_{j - 1}\}$ is adjacent to a vertex in $\{x_{j+1},\dots,x_{\ell}\}$ (by Lemma~\ref{diagonallemma} as $x_1x_j$ is a diagonal of $G_1$), the claim follows.\hfill{\smallqed}

\parskip=3mm
Let $G'$ be the graph obtained from $G$ by deleting the vertices $x_2,x_3,\dots,x_{\ell}$. Then, $G'$ is the mop $G_2$. Let $n' = |V(G')|$. Then, $n' = n - (\ell - 1) \le n - (k + 4)$. Suppose $n' \le k+3$. Then, $G_2 - \{x_1, x_{\ell + 1}\}$ contains no copy of $K_{1,k+1}$. Together with $x_1, x_{\ell + 1} \in N_{G}[x_j]$ and Claim~\ref{c2}, this gives us that $\{x_j\}$ is a $K_{1, k + 1}$-isolating set of $G$, so $\ii_k(G) \le 1 < n/(k + 4)$. Now suppose $n' \ge k+4$. By the induction hypothesis, $\ii_k(G') \le n'/(k + 4) \le n/(k + 4)-1$. Let $S'$ be a $K_{1, k + 1}$-isolating set of $G'$ with $|S'| = \ii_k(G')$. Since $x_1, x_{\ell + 1} \in N_{G}[x_j]$, it follows by Claim \ref{c2} that $S' \cup \{x_{j}\}$ is a $K_{1, k + 1}$-isolating set of $G$, so $\ii_k(G) \le |S'| + 1 \le n/(k + 4)$.~\hfill{\qed}

\parskip=5mm



\noindent
\textbf{Proof of Theorem~\ref{t:main2}.} We use an inductive argument similar to that in the proof of Theorem~\ref{t:main1}. 

\parskip=4pt
If $k+3 \le n \le 2k+7$, then, by Lemmas~\ref{mop:iotaofsmallmop} and \ref{basicfacts}(b), $\ii_k(G) \leq 1 \le (n+n_2)/(k+5)$. We now assume that $n \ge 2k+8$ and proceed by induction on $n$. Let $x_1x_2 \ldots x_nx_1$ be the unique Hamiltonian cycle $C$ of $G$ and hence the boundary of the exterior face of $G$. By Lemma~\ref{mop:6-10} with $r = k+2$, $G$ has a diagonal $d$ that partitions it into two mops $G_1$ and $G_2$ such that $G_1$ has exactly $\ell$ Hamiltonian edges of $G$ for some $\ell \in [k+4, 2k+6]$. We may assume that $d = x_1x_{\ell + 1}$ and $V(G_1) =  \{x_1,x_2,\ldots,x_{\ell + 1}\}$. Note that $x_1x_2,x_2x_3,\dots,x_{\ell}x_{\ell + 1}$ are the $\ell$ Hamiltonian edges of $G$ that belong to $G_1$. Let $(x_1,x_j,x_{\ell + 1})$ be the triangular face of $G_1$ that contains the edge $x_1x_{\ell + 1}$. Then, $2 \le j \le \ell$.

Suppose $\ell = k + 4$. Then, $x_{k+5},x_n \in N_{G_2}(x_1)$ and $x_1,x_{k+6} \in N_{G_2}(x_{k+5})$. Since $x_1$ and $x_{k+5}$ are adjacent in $G_2$, Lemma~\ref{maxofn2} tells that the degrees of $x_1$ and $x_{k+5}$ in $G_2$ cannot both be $2$. Thus, $d_{G_2}(x_1) + d_{G_2}(x_{k + 5}) \ge 5$.

Suppose $d_{G_2}(x_1) + d_{G_2}(x_{k + 5}) = 5$. We may assume that $d_{G_2}(x_1) = 3$ and $d_{G_2}(x_{k + 5}) = 2$.  We have $x_1x_n, x_{k+5}x_{k+6} \in E(C) \cap E(G_2)$. Since $x_1x_{k+5}, x_{k+5}x_{k+6} \in E(G_2)$ and $d_{G_2}(x_{k+5}) = 2$, $N_{G_2}(x_{k+5}) = \{x_1,x_{k+6}\}$. Thus, since $G_2$ is a mop, the face having $x_1x_{k+5}$ and $x_{k+5}x_{k+6}$ on its boundary must also have $x_1x_{k+6}$ on its boundary (as all interior faces are triangles), that is, $x_1x_{k+6} \in E(G_2)$ (see Figure 4). Together with $x_1x_{k+5}, x_1x_n \in E(G_2)$ and $d_{G_2}(x_1) = 3$, this gives us $N_{G_2}(x_1) = \{x_{k+5},x_{k+6},x_n\}$. Thus, since $G_2$ is a mop, the face having $x_1x_{k+6}$ and $x_1x_n$ on its boundary must also have $x_{k+6}x_n$ on its boundary, that is, $x_{k+6}x_n \in E(G_2)$. Let $G' = G - \{x_1,x_2,\dots,x_{k+5}\}$. Then, $G' = G_2 - \{x_1, x_{k+5}\}$. Since $d_{G_2}(x_{k+5}) = 2$, $G_2 - x_{k+5}$ is a mop by Lemma~\ref{basicfacts}(c). Since $d_{G_2 - x_{k+5}}(x_1) = 2$, $G'$ is a mop by Lemma~\ref{basicfacts}(c). Let $n' = |V(G')|$ and $n_2' = |\{v \in V(G') \colon d_{G'}(v) = 2\}|$. We have $n' = n - (k+5) \ge k+3$.
\vskip 5 pt

\setlength{\unitlength}{0.8cm}
\begin{center}
\begin{picture}(2, 6)
\put(-1, 0){\circle*{0.2}}
\put(3, 0){\circle*{0.2}}
\put(-2, 2){\circle*{0.2}}
\put(-1, 4){\circle*{0.2}}
\put(3, 4){\circle*{0.2}}
\put(4, 2){\circle*{0.2}}
\put(1, 5){\circle*{0.2}}
\put(-1, 0){\line(-1, 2){1}}
\put(-1, 4){\line(-1, -2){1}}
\put(-1, 4){\line(2, 1){2}}

\put(4, 2){\line(-1, 2){1}}
\put(4, 2){\line(-1, -2){1}}
\put(-1.8, 0){\small{$x_1$}}
\put(-2.8, 2){\small{$x_2$}}
\put(-1.4, 4.2){\small{$x_3$}}
\put(1, 5.2){\small{$x_4$}}
\put(1.8, 4.4){\small{$\cdot$}}
\put(2, 4.3){\small{$\cdot$}}
\put(2.2, 4.2){\small{$\cdot$}}
\put(3, 4.2){\small{$x_{k + 3}$}}
\put(4.2, 2){\small{$x_{k + 4}$}}
\put(3.3, 0){\small{$x_{k + 5}$}}
\put(3.3, -1.7){\small{$x_{k + 6}$}}
\put(-1.8, -1.7){\small{$x_n$}}

\put(-1, -2){\circle*{0.2}}
\put(3, -2){\circle*{0.2}}
\put(3, 0){\circle*{0.2}}
\put(3, 0){\circle*{0.2}}

\put(-1, 0){\line(1, 0){4}}
\put(-1, -2){\line(0, 1){2}}
\put(-3, -2){\line(1, 0){8}}
\put(3, -2){\line(0, 1){2}}
\put(3, -2){\line(-2, 1){4}}
\end{picture}
\end{center}
\begin{center} \vskip 58 pt \footnotesize{\textbf{\textbf{Figure 4}}}
\end{center}
\vskip 30 pt

By Lemma~\ref{maxofn2}, at most one of $x_{k+6}$ and $x_n$ has degree~$2$ in $G'$. By Lemma~\ref{maxofn2}, at most one of $x_1$ and $x_{k+5}$ has degree~$2$ in $G_1$, and hence, by Lemma~\ref{basicfacts}(b), $d_{G_1}(x_h) = 2$ for some $h \in [2, k+4]$. Since $x_h \in V(G_1) \setminus V(G_2)$, $d_{G}(x_h) = d_{G_1}(x_h)$. Therefore, $n_2' \le n_2$, and hence $n'+n_2' \le n+n_2-(k+5)$. By the induction hypothesis, $\ii_k(G') \le (n'+n_2')/(k+5) \le (n+n_2)/(k+5) - 1$. Let $S'$ be a smallest $K_{1,k+1}$-isolating set of $G'$. Clearly, $|V(G_1 - N_{G_1}[x_j])| \le k+1$, so $G_1 - N_{G_1}[x_j]$ does not contain a copy of $K_{1, k+1}$. Since $x_j$ is adjacent to both $x_1$ and $x_{k+5}$, it follows that $S' \cup \{x_j\}$ is a $K_{1,k+1}$-isolating set of $G$. Thus, we have $\ii_k(G) \le |S'| + 1 = \ii_k(G') + 1 \le (n+n_2)/(k+5)$.

Now suppose $d_{G_2}(x_1) + d_{G_2}(x_{k+5}) \ge 6$. Let $G'$ be the graph obtained from $G$ by deleting the vertices $x_2,x_3,\dots,x_{k+4}$ and contracting the edge $x_1x_{k+5}$ to form a new vertex $y$. Then, $G'$ is obtained from $G_2$ by contracting $x_1x_{k+5}$. Thus, $G'$ is a mop by Lemma~\ref{contraction}. Let $n' = |V(G')|$ and $n_2' = |\{v \in V(G') \colon d_{G'}(v) = 2\}|$. We have $n' = n - (k+4) \ge k+4$.

Suppose $d_{G'}(y) \le 2$. As noted above, $x_1x_n, x_{k+5}x_{k+6} \in E(G_2)$. Thus, $N_{G'}(y) = \{x_{k+6}, x_n\}$. Since $d_{G_2}(x_1) + d_{G_2}(x_{k+5}) \ge 6$, we obtain $N_{G_2}(x_1) = \{x_{k+5}, x_{k+6}, x_n\}$ and $N_{G_2}(x_{k+5}) = \{x_1, x_{k+6}, x_n\}$. Since $N_{G_2}(x_1) = \{x_{k+5}, x_{k+6}, x_n\}$, $x_1x_{k+6}$ is a diagonal of $G_2$. By Lemma~\ref{diagonallemma}, we obtain $x_{k+5}x_n \notin E(G_2)$, which contradicts $N_{G_2}(x_{k+5}) = \{x_1, x_{k+6}, x_n\}$. Therefore, $d_{G'}(y) \ge 3$.

Suppose that every vertex that has degree~$2$ in $G'$ also has degree~$2$ in $G$. As in the proof for the case $d_{G_2}(x_1) + d_{G_2}(x_{k + 5}) = 5$, $d_{G_1}(x_h) = 2$ for some $h \in [2, k+4]$, so $n_2' \le n_2 - 1$. Thus, $n'+n_2' \le n+n_2-(k+5)$. By the induction hypothesis, $\ii_k(G') \le (n'+n_2')/(k+5) \le (n+n_2)/(k+5) - 1$. Let $S'$ be a smallest $K_{1,k+1}$-isolating set of $G'$. Then, $|S'| \le (n+n_2)/(k+5) - 1$. We can continue as in the proof of Theorem~\ref{t:main1} for the case $\ell = k+4$ to obtain $\ii_k(G) \le |S'| + 1 \le (n+n_2)/(k+5)$.

Now suppose that $G'$ has a vertex $z$ such that $d_{G'}(z) = 2 \neq d_{G}(z)$. As $z \neq y$, we have $x_{1}, x_{k+5} \in N_{G_{2}}(z)$. For each $i \in [k+7, n-1]$ with $d_{G'}(x_i) = 2$, we have $N_{G'}(x_i) = N_G(x_i) = \{x_{i-1}, x_{i+1}\}$, so $z \neq x_i$. Thus, $z = x_{k+6}$ or $z = x_n$. By symmetry, we may assume that $z = x_{k+6}$. Since $x_{k+6}x_{k+7} \in E(C) \cap E(G_2)$ and $x_1, x_{k+5} \in N_{G_2}(x_{k+6})$, $N_{G_2}(x_{k+6}) = \{x_1, x_{k+5}, x_{k+7}\}$. Thus, since $G_2$ is a mop, the face having $x_1x_{k+6}$ and $x_{k+6}x_{k+7}$ on its boundary must also have $x_1x_{k+7}$ on its boundary (as all interior faces are triangles), meaning that $x_1x_{k+7} \in E(G_2)$. By Lemma~\ref{diagonallemma}, $x_1x_{k+7}$ partitions $G$ into two mops $H_1$ and $H_2$ such that $V(H_2) = \{x_1, x_{k+7}, x_{k+8}, \dots, x_n\}$. Let $G^* = H_2$, $n^* = |V(H_2)|$, and $n_2^* = |\{v \in V(H_2) \colon d_{H_2}(v) = 2\}|$. We have $n^* = n-(k+5) \geq k+3$. By Lemma~\ref{maxofn2}, for each $i \in \{1, 2\}$, at most one of $x_1$ and $x_{k+7}$ has degree~$2$ in $H_i$. By Lemma~\ref{basicfacts}(b), $d_{H_1}(x_h) = 2$ for some $h \in V(H_1) \setminus \{x_1, x_{k+7}\}$, and hence $d_{G}(x_h) = 2$. Therefore, $n_2^* \le n_2$, and hence $n^*+n_2^* \le n+n_2-(k+5)$. By the induction hypothesis, $\ii_k(G^*) \le (n^*+n_2^*)/(k+5) \le (n+n_2)/(k+5) - 1$. Let $S^*$ be a smallest $K_{1,{k+1}}$-isolating set of $G^*$. Let $x^* = x_{k+5}$ if $j = 2$, and let $x^* = x_1$ otherwise. If $3 \leq j \leq k+4$, then $x_1, x_2, x_j, x_{k+5}, x_{k+6}, x_{k+7} \in N_{H_1}[x^*]$. If $j = 2$, then $x_1, x_2, x_{k+4}, x_{k+5}, x_{k+6} \in N_{H_1}[x^*]$. Since $x_1x_{k+5}$ is a diagonal of $G$, we have $\{x_2, x_3, \dots, x_{k+4}\} \cap N_G(x_{k+7}) = \emptyset$ by Lemma~\ref{diagonallemma}. Therefore, $S^* \cup \{x^*\}$ is a $K_{1,{k+1}}$-isolating set of $G$. Thus, we have $\ii_k(G) \le |S^*| + 1 = \ii_k(G^*) + 1 \le (n+n_2)/(k+5)$.

Having settled the case $\ell = k+4$, we now use induction on $\ell$. Thus, we consider $\ell \ge k + 5$ and assume that, if $G$ has a diagonal that partitions it into two mops $H_1$ and $H_2$ such that $H_1$ has exactly $\ell^*$ Hamiltonian edges of $G$ for some $\ell^* \in [k+4, \ell - 1]$, then $\ii_k(G) \leq (n + n_2)/(k+5)$. Since $\ell \le 2k+6$, $\ell = k + 4 + t$ for some $t \in [k + 2]$. By the argument in the proof of Claim~\ref{c1}, if $j \notin [t + 2, k + 4]$, then $\ii_k(G) \leq (n+n_2)/(k+5)$. Now suppose $j \in [t + 2, k + 4]$. By the same argument for Claim~\ref{c2}, $G_1 - N_{G_{1}}[x_j]$ contains no copy of $K_{1,k+1}$.

Let $G'$ be the graph obtained from $G$ by deleting the vertices $x_2,x_3,\dots,x_{\ell}$. Then, $G'$ is the mop $G_2$. Let $n' = |V(G')|$ and $n_2' = |\{v \in V(G') \colon d_{G'}(v) = 2\}|$. Then, $n' = n - (\ell - 1) \le n - (k + 4)$. Suppose $n' \le k+3$. Since $x_1, x_{\ell + 1} \in N_G(x_j) \cap V(G')$ and $G_1 - N_{G_{1}}[x_j]$ contains no copy of $K_{1,k+1}$, $\{x_j\}$ is a $K_{1, k + 1}$-isolating set of $G$, so $\ii_k(G) \le 1 < (n+n_2)/(k + 5)$. Now suppose $n' \ge k+4$. Since $x_1x_{\ell + 1}$ is a diagonal of $G$, Lemma~\ref{diagonallemma} gives us that $d_{G'}(v) = d_G(v)$ for each $v \in V(G') \setminus \{x_1, x_{\ell + 1}\}$. By Lemma~\ref{maxofn2}, at most one of $x_1$ and $x_{\ell + 1}$ has degree~$2$ in $G'$. We have $x_1, x_j, x_{\ell} \in N_{G_1}(x_{\ell + 1})$. Since $x_1x_j$ is a diagonal of $G_1$ (as $t+2 \le j \le k+4 < \ell$), Lemmas~\ref{diagonallemma}, \ref{maxofn2}, and \ref{basicfacts}(b) give us that at least one vertex in $\{x_2,\dots,x_{j-1}\}$ has degree~$2$ in $G$, and that at least one vertex in $\{x_{j+1},\dots,x_{\ell}\}$ has degree~$2$ in $G$. Thus, $n'_2 \le n_2-1$, and hence $n'+n'_2 \le n+n_2 - (k+5)$. Let $S'$ be a smallest $K_{1,k+1}$-isolating set of $G'$. By the induction hypothesis, $|S'| \le (n'+n'_2)/(k + 5) \le (n+n_2)/(k + 5)-1$. Since $x_1, x_{\ell + 1} \in N_G(x_j)$ and $G_1 - N_{G_{1}}[x_j]$ contains no copy of $K_{1,k+1}$,  $S' \cup \{x_{j}\}$ is a $K_{1,k + 1}$-isolating set of $G$, so $\ii_k(G) \le |S'| + 1 \le (n+n_2)/(k + 5)$.~\hfill{\qed}

\parskip=5mm

\noindent
\textbf{Proof of Theorem~\ref{t:main3}.} By (\ref{n_2bounds}), $2 \leq n_2 \le n/2$. Since $n - n_2 \ge n/2 \ge (2k+3)/2$, $n - n_2 \geq k+2$. Let $V_2$ be the set of vertices of $G$ of degree~$2$, let $G' = G - V_2$, and let $n' = |V(G')|$. Then, $n' = n - n_2$. By Lemmas~\ref{maxofn2} and \ref{basicfacts}(c), $G'$ is a mop. If $k = 1$, then let $S$ be a smallest dominating set of $G'$. If $k \geq 2$, then let $S$ be a smallest $K_{1,k-1}$-isolating set of $G'$. By Theorem~\ref{t:known1}(a) and Theorem~\ref{t:main1}, $|S| \le n'/(k+2)$. By Lemmas~\ref{maxofn2} and \ref{basicfacts}(e), $V_2$ is an independent set of $G$ and, in $G$, each vertex in $V(G')$ is adjacent to at most two vertices in $V_2$. Consequently, $S$ is a $K_{1,k+1}$-isolating set of $G$, and hence $\ii_k(G) \le |S| \le (n - n_2)/(k+2)$.~\hfill{\qed}

\parskip=5mm

\noindent
\textbf{Proof of Theorem~\ref{t:main4}.} 
If $n = 4$, then the result is trivial. Suppose $n \ge 5$. Let $V_2$ be the set of vertices of $G$ of degree~$2$, let $G' = G- V_2$, and let $n' = |V(G')|$. Then, $n' = n - n_2$. By Lemma~\ref{maxofn2}, $n \ge 2n_2$, so $n' \ge n/2$. Since $n \ge 5$, $n' \geq 3$. By Lemmas~\ref{maxofn2} and \ref{basicfacts}(c), $G'$ is a mop. Let $x_0x_1 \dots x_{n'-1}x_0$ be the unique Hamiltonian cycle $C'$ of $G'$. Let $y_1, \dots, y_{n_2}$ be the vertices in $V_2$. 

\parskip=0pt

\begin{lem}
\label{V_2 lemma}
$N_G(y_1), \dots, N_G(y_{n_2})$ are distinct edges of $C'$.
\end{lem}
\proof Consider any $r \in [n_2]$. By Lemma~\ref{maxofn2}, $N_G(y_r) \subseteq V(G')$. Since $y_r \in V_2$, $N_G(y_r) = \{x_h, x_{(h + j) \bmod n'}\}$ for some $h, j \in \{0, 1, \dots, n'-1\}$ with $j \geq 1$ (where $\bmod$ is the usual modulo operation). Since $G$ is a mop, no vertex of $G$ lies in the interior of a cycle of $G$. Thus, none of $y_1, \dots, y_{n_2}$ lie in the interior of $C'$. 
Suppose $2 \le j \le n'-2$. Since $x_{(h+n'-1) \bmod n'}$ does not lie in the interior of the cycle 
\[x_h x_{(h+1) \bmod n'} \dots x_{(h+j) \bmod n'} y_r x_h\]
of $G$, we obtain that $x_{(h+1) \bmod n'}$ lies in the interior of the cycle 
\[x_h y_r x_{(h+j) \bmod n'} x_{(h+j+1) \bmod n'} \dots x_{(h+n'-1) \bmod n'}x_h\]
of $G$, a contradiction. Thus, $j \in \{1, n'-1\}$, and hence $N_G(y_r) \in E(C')$. We may assume that $j = 1$. Suppose $N_G(y_s) = N_G(y_r)$ for some $s \in [n_2] \setminus \{r\}$. Since $x_{(h+2) \bmod n'}, \dots, x_{(h+n'-1) \bmod n'}$ do not lie in the interior of the cycle $x_h x_{(h+1) \bmod n'} y_r x_h$ of $G$ and do not lie in the interior of the cycle $x_h x_{(h+1) \bmod n'} y_s x_h$ of $G$, and $y_r$ does not lie in the interior of the cycle
%
%
$x_h y_s x_{(h+1) \bmod n'}x_h$ of $G$, we obtain that $y_s$ lies in the interior of the cycle 
%
%
$x_h y_r x_{(h+1) \bmod n'}x_h$ of $G$, a contradiction. Thus, $N_G(y_s) \neq N_G(y_r)$. Therefore, the lemma is proved. \hfill{\smallqed}

\parskip=4pt

Suppose $e' \notin \{N_G(y_1), \dots, N_G(y_{n_2})\}$ for some edge $e'$ of $C'$. By Lemma~\ref{V_2 lemma}, $n_2 < n'$, so $n > 2n_2$. We may assume that $e' = x_1x_2$. Then, $\{x_0\} \cup \{x_{2i+1} \colon 1 \le i \le \lfloor n'/2 \rfloor - 1\}$ is a dominating set of $G$, and hence $\gamma(G) \leq n'/2$.

Now suppose $e \in \{N_G(y_1), \dots, N_G(y_{n_2})\}$ for each edge $e$ of $C'$. By Lemma~\ref{V_2 lemma}, $n_2 = n'$, so $n = 2n_2$. If $n_2$ is even, then $\{x_{2i-1} \colon 1 \le i \le n'/2\}$ is a dominating set of $G$, so $\gamma(G) \leq n'/2$. If $n_2$ is odd, then $\{x_0\} \cup \{x_{2i-1} \colon 1 \le i \le (n'-1)/2\}$ is a dominating set of $G$, so $\gamma(G) \leq (n' + 1)/2$.

We have verified the bound in the theorem. We now show that the bound is attainable. 
%

Suppose that $G$ is as in Definition~\ref{extreme}. 
The Hamiltonian cycle of $G$ is $x_1y_1x_2y_2 \dots x_ty_tx_1$. We have $V_2 = \{y_i \colon i \in [t]\}$, $n_2 = t$, and $n = 2n_2$. Let $s = \lfloor (t-1)/2 \rfloor$. Let $S$ be a smallest dominating set of $G$. Let $D = (S \setminus \{y_i \colon i \in [t]\}) \cup \{x_i \colon i \in [t], \, y_i \in S\}$. Then, $D \subseteq \{x_i \colon i \in [t]\}$, $|D| \le |S|$, and $D$ is a dominating set of $G$. We may assume that $x_1 \in D$. For each $i \in [t-1]$, we have $y_i \in N_G[D]$, so at least one of $x_i$ and $x_{i+1}$ is in $D$. Thus, since $D \supseteq \{x_1\} \cup \bigcup_{i=1}^s (D \cap \{x_{2i}, x_{2i+1}\})$, we obtain $|D| \ge 1 + \sum_{i=1}^s |D \cap \{x_{2i}, x_{2i+1}\}| \geq 1 + s$. Since $\{x_{2i} \colon i \in [s]\} \cup \{x_t\}$ is a dominating set of $G$, $\gamma(G) \le 1 + s$. Thus, since $1 + s \leq |D| \leq |S| = \gamma(G)$, $\gamma(G) = 1+s$. We have $s = \lfloor (n_2 - 1)/2 \rfloor$. If $n_2$ is odd, then $\gamma(G) = 1 + (n_2 - 1)/2 = (n - n_2 + 1)/2 = (n - n_2 + {\bf 1}(n,n_2))/2$. If $n_2$ is even, then $\gamma(G) = 1 + (n_2 - 2)/2 = (n - n_2)/2 = (n - n_2 + {\bf 1}(n,n_2))/2$.~\hfill{\qed}


\parskip=5mm

\noindent
\textbf{Proof of Theorem~\ref{t:MAIN}.} The bound and part (a) are immediate consequences of Theorems~\ref{t:known2.0} and~\ref{t:main4}. 

\parskip=4pt

Let $V_2$ be the set of vertices of $G$ of degree~$2$, let $y_1, \dots, y_{n_2}$ be the members of $V_2$, let $G' = G - V_2$, and let $n' = |V(G')|$. Then, $n' = n - n_2$. By Lemma~\ref{basicfacts}(b), $n_2 \geq 2$. If $n \ge 5$, then $n' \geq 3$ by Lemma~\ref{maxofn2}, so $G'$ is a mop by Lemmas~\ref{maxofn2} and \ref{basicfacts}(c). Let $C'$ be the Hamiltonian cycle of $G'$. By Lemma~\ref{V_2 lemma}, $N_G(y_1), \dots, N_G(y_{n_2})$ are distinct edges of $C'$ if $n \ge 5$. For each $x \in V(G')$ with $d_{G'}(x) = 2$, $x \in N_G[V_2]$ as $x \notin V_2$.

Suppose $\gamma(G) = n/3$. By the first part of the result, either $n_2 = n/3$ or $(n,n_2) = (6,3)$. Thus, $n \geq 6$. If $n_2 = n/3$, then $n' = 2n_2$, so $G$ is $n_2$-special. If $(n,n_2) = (6,3)$, then $C'$ is a $3$-vertex cycle and $N_G(y_1), \dots, N_G(y_{n_2})$ are its $3$ edges, so $G$ is $3$-extreme. Thus, (b) is proved.

Suppose that $G$ is extra $n_2$-special. Then, $n = \sum_{i=1}^{n_2} |N_G[y_i]| = 3n_2$. Thus, $n \geq 6$. If $S$ a dominating set of $G$, then $S \cap N_G[y_i] \neq \emptyset$ for each $i \in [n_2]$. Thus, $\gamma(G) \ge n_2$. Since $V_2$ is a dominating set of $G$, $\gamma(G) = n_2$. Since $n = 3n_2$, (c) follows.~\hfill{\qed}

\section{An Art Gallery Theorem relaxation for guarding corners} \label{AGTextsection}
For an integer $k \geq -1$ and a mop $G$, let $\mathcal{P}_k(G)$ denote the set of $(k+2)$-vertex paths of the Hamiltonian cycle of $G$, and let $c_{k}(G)$ denote the size of a smallest set $S$ such that no member of $\mathcal{P}_k(G)$ is a subgraph of $G - N[S]$. By the same argument used in the proof of Theorem~\ref{t:main1} (including that of Lemma~\ref{mop:iotaofsmallmop}), we have the following theorem.

\begin{thm}\label{thmconsecutive} If $k \ge -1$, $n \ge k + 4$, and $G$ is a mop of order $n$, then \[c_{k}(G) \le \frac{n}{k + 4}.\]
\end{thm}

\noindent It is worth pointing out that the same argument yields stronger results; in particular, we have that if $\mathcal{F}_k$ is the set of connected $(k+2)$-vertex graphs, then $\iota(G, \mathcal{F}_k) \le n/(k+4)$. 

Theorem~\ref{thmconsecutive} gives a result in computational geometry that extends the Art Gallery Theorem. We assume that an art gallery is the closed set of points bounded by a polygon $P$ of $n$ sides (so $P$ has $n$ corners). Two points in $P$~(including the sides and corners of $P$ as the set is closed) are \emph{visible} if the straight line joining them does not intersect the exterior of $P$. The classical problem solved by Chv\'{a}tal~\cite{Ch75} was to find the minimum number of guards that can be placed in $P$ so that every point in $P$ is visible to at least one guard. We relax the problem by restricting the visibility condition to corners only and allowing the guards to ignore sets of at most $k + 1$ consecutive corners on the perimeter of $P$ in the following sense: more than $k + 1$ corners may be ignored, but at least one of every $k + 2$ consecutive corners needs to be visible to at least one guard. Let $g_{k}(P)$ denote the minimum number of guards that can be used for this purpose. Note that having $k = -1$ means that we do not allow the guards to ignore any corner. With a slight abuse of notation, let $V(P)$ denote the set of corners of $P$. Using Theorem~\ref{thmconsecutive}, we obtain the following result.

\begin{thm} \label{AGTextthm}
If $k \ge -1$, $n \ge k + 4$, and $n$ is the number of corners of a polygon $P$, then
%
\[g_{k}(P) \le \frac{n}{k + 4}.\]
Moreover, for every $t \ge 1$, there exists a polygon $P_{k,t}$ such that $g_{k}(P_{k,t}) = t = \frac{|V(P_{k,t})|}{k + 4}$.
\end{thm}

\proof We may represent $P$ by a cycle $C_P$ drawn on the plane. Thus, the vertices of $C_P$ represent the corners of $P$, and the edges of $C_P$ represent the sides of $P$. We insert non-crossing edges in the interior of $C_P$ (without adding vertices) until we obtain a mop $G_P$. Thus, $C_P$ is the Hamiltonian cycle of $G_P$. 
By Theorem~\ref{thmconsecutive}, there exists a set $S$ of at most $n/(k + 4)$ vertices of $G_P$ such that $G_P - N_{G_P}[S]$ contains no $(k + 2)$-vertex path of $C_P$. By placing guards at the corners of $P$ represented by the vertices in $S$, we obtain that at least one of every $k + 2$ consecutive corners of $P$ is visible to at least one guard. Therefore, $g_{k}(P) \le n/(k+4)$.

We now show that the bound in the theorem can be attained for any integer value of $n/(k+4)$ by constructing $P_{k,t}$ explicitly. If the interior angle at a corner of a polygon is reflex (that is, more than $180^{\circ}$), then we call the corner a \emph{reflex corner}. A \emph{reflex chain} is a sequence of consecutive reflex corners. A polygon is \emph{spiral} if it is a triangle or has exactly one reflex chain. Let $t$ be a positive integer. For each $i \in [t]$, let $S_{i}$ be a spiral polygon of $k + 4$ sides, and let $c^{i}_{1}, \dots, c^{i}_{k + 4}$ be the corners of $S_{i}$, listed in the order they appear in the clockwise direction and such that, if $k \ge 0$, then $c^{i}_{2}, \dots, c^{i}_{k}$ is the reflex chain. If $t = 1$, then we take $P_{k,t}$ to be $S_1$, and we trivially have $g_{k}(P_{k, t}) = 1 = |V(P_{k,t})|/(k+ 4)$. Suppose $t \ge 2$. Place $S_{1}, \dots, S_{t}$ consecutively on a plane in such a way that no two intersect, $c^{1}_{k + 4}, c^{2}_{1}, c^{2}_{k + 4}, \dots, c^{t - 1}_{1}, c^{t - 1}_{k + 4}$ and $c^{t}_{1}$ are on the same horizontal line $L_{1}$, and $c^{1}_{1}$ and $c^{t}_{k + 4}$ are on the same horizontal line $L_{2}$ slightly below $L_{1}$. For each $i \in [t]$, remove the side $c^{i}_{1}c^{i}_{k + 4}$ of $S_{i}$. Join $c^{1}_{1}$ to $c^{t}_{k + 4}$ by a line segment, and for each $i \in [t - 1]$, join $c^{i}_{k + 4}$ to $c^{i + 1}_{1}$ by a line segment. Let $P_{k, t}$ be the polygon obtained. The polygon $P_{3, 4}$ is illustrated in Figure~$5$. 
We have $|V(P_{k,t})| = t(k+4)$. For each $i \in [t]$, if we place a guard at the corner $c^{i}_{1}$ of $P_{k,t}$, then each corner of $S_i$ that is not visible to the guard is one of the $k + 1$ consecutive corners $c^{i}_{3}, \dots, c^{i}_{k + 3}$. Consequently, 
$g_{k}(P_{k, t}) \le t$. We take $L_{1}$ and $L_{2}$ close enough so that each of the sets $V(S_{1}), \dots, V(S_{t})$ of consecutive corners of $P_{k,t}$ needs its own guard. 
Therefore, $g_{k}(P_{k, t}) = t$.~\hfill{\qed}

\vskip 5 pt

\begin{figure}[h]
\centering
\includegraphics[width=12cm]{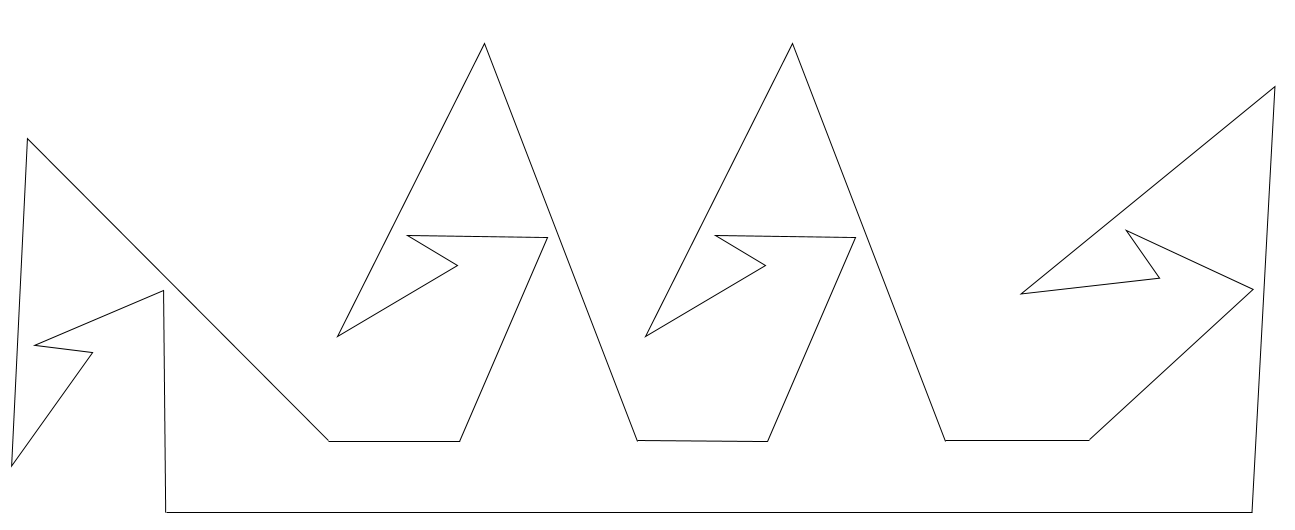}
\end{figure}
\begin{center}
\footnotesize{\textbf{Figure 5 :} The polygon $P_{3, 4}$.}
\end{center}
\vskip 30 pt

\noindent \textbf{Acknowledgements.} 
The authors are grateful to two anonymous referees for checking the paper and providing constructive remarks. The second author was supported in part by Skill Development Grant 2018, King Mongkut's University of Technology Thonburi.


\begin{thebibliography}{99}

\bibitem{AZ} M. Aigner and G. Ziegler, Proofs from THE BOOK (6th ed.), Springer, Berlin/Heidelberg, 2018.

\bibitem{Borg_conn} P. Borg. Isolation of connected graphs, arXiv:2110.03773 [math.CO].

\bibitem{Borg1} P. Borg, Isolation of cycles, \textit{Graphs Combin.} \textbf{36} (2020), 631--637.

\bibitem{BFK} P. Borg, K. Fenech and P. Kaemawichanurat, Isolation of $k$-cliques, \textit{Discrete Math.} \textbf{343} (2020), paper 111879.

\bibitem{BFK2} P. Borg, K. Fenech and P. Kaemawichanurat, Isolation of $k$-cliques II, \textit{Discrete Math.}, in press.

\bibitem{BK} P. Borg and P. Kaemawichanurat, Partial domination of maximal outerplanar graphs, \textit{Discrete Appl. Math.} \textbf{283} (2020), 306--314.

\bibitem{CaHa17} Y. Caro and A. Hansberg, Partial domination - the isolation number of a graph, \textit{FiloMath} \textbf{31:12} (2017), 3925--3944.

\bibitem{CaWa13} C.N. Campos and Y. Wakabayashi, On dominating sets of maximal outerplanar graphs, \textit{Discrete Appl. Math.} \textbf{161} (2013), 330--335.

\bibitem{CaCaHeMa16} S. Canales, I. Castro, G. Hern\'{a}ndez and M. Martins, Combinatorial bounds on connectivity for dominating sets in maximal outerplanar graphs, \textit{Electron. Notes Discrete Math.} \textbf{54} (2016), 109--114.

\bibitem{Ch75} V. Chv\'{a}tal, A combinatorial theorem in plane geometry, \textit{J. Combin. Theory Ser. B} \textbf{18} (1975), 39--41.

\bibitem{C} E.J. Cockayne, Domination of undirected graphs -- A survey, Lecture Notes in Mathematics, Volume 642, Springer, 1978, 141--147.

\bibitem{CH} E.J. Cockayne and S.T. Hedetniemi, Towards a theory of domination in graphs, \textit{Networks} \textbf{7} (1977), 247--261.

\bibitem{DoHaJo16} M. Dorfling, J.H. Hattingh and E. Jonck, Total domination in maximal outerplanar graphs II, \textit{Discrete Math.} \textbf{339} (2016), 1180--1188.

\bibitem{DoHaJo17} M. Dorfling, J.H. Hattingh and E. Jonck, Total domination in maximal outerplanar graphs, \textit{Discrete Appl. Math.} \textbf{217} (2017), 506--511.

\bibitem{FK} O. Favaron and P. Kaemawichanurat, Inequalities between the $K_k$-isolation number and the independent $K_k$-isolation number of a graph, \textit{Discrete Appl. Math.} \textbf{289} (2021), 93--97.

\bibitem{Fi78} S. Fisk, A short proof of Chv\'{a}tal's watchman theorem, \textit{J. Combin. Theory Ser. B} \textbf{24} (1978), 374.

\bibitem{HHS} T.W. Haynes, S.T. Hedetniemi and P.J. Slater, Fundamentals of Domination in Graphs, Marcel Dekker, Inc., New York, 1998.

\bibitem{HHS2} T.W. Haynes, S.T. Hedetniemi and P.J. Slater (Editors), Domination in Graphs: Advanced Topics, Marcel Dekker, Inc., New York, 1998.

\bibitem{HL} S.T. Hedetniemi and R.C. Laskar (Editors), Topics on Domination, \textit{Discrete Math.} \textbf{86} (1990).

\bibitem{HL2} S.T. Hedetniemi and R.C. Laskar, Bibliography on domination in graphs and some basic definitions of domination parameters, \textit{Discrete Math.} \textbf{86} (1990), 257--277.

\bibitem{HeKa18} M. A. Henning and P. Kaemawichanurat, Semipaired domination in maximal outerplanar graphs, \textit{J. Comb. Optim.} \textbf{38} (2019), 911--926.





\bibitem{LeZuZy17} M. Lema\'{n}ska, R. Zuazua and P. Zylinski, Total dominating sets in maximal outerplanar graphs, \textit{Graphs Combin.} \textbf{33} (2017), 991--998.

\bibitem{Li16} Z. Li, E. Zhu, Z. Shao and J. Xu, On dominating sets of maximal outerplanar and planar graphs, \textit{Discrete Appl. Math.} \textbf{198} (2016), 164--169.

\bibitem{MaTa96} L. R. Matheson and R. E. Tarjan, Dominating sets in planar graphs, \textit{European J. Combin.} \textbf{17} (1996), 565--568.

\bibitem{Ro83} J. O'Rourke, Galleries need fewer mobile guards: a variation to Chv\'{a}tal's theorem, \textit{Geom. Dedicata} \textbf{14} (1983), 273--283.

\bibitem{Ro87} J. O'Rourke, Art gallery theorems and algorithms, Oxford University Press, New York, 1987.

\bibitem{To13} S. Tokunaga, Dominating sets of maximal outerplanar graphs, \textit{Discrete Appl. Math.} \textbf{161} (2013), 3097--3099.

\bibitem{KaJi} S. Tokunaga, T. Jiarasuksakun and P. Kaemawichanurat, Isolation number of maximal outerplanar graphs, \textit{Discrete Appl. Math.} \textbf{267} (2019), 215--218.

\bibitem{Y} J. Yan, Isolation of the diamond graph, \textit{Bull. Malaysian Math. Sci. Soc.} (2022), in press.

\bibitem{YW} H. Yu and B. Wu, Admissible property of graphs in terms of radius, \textit{Graphs Combin.} \textbf{38} (2022), paper 6.

\bibitem{ZW} G. Zhang and B. Wu, $K_{1,2}$-isolation in graphs, \textit{Discrete Appl. Math.} \textbf{304} (2021), 365--374.

\end{thebibliography}
\end{document}